\documentclass[final]{siamart1116}

\usepackage[T1]{fontenc}
\usepackage{lipsum}
\usepackage{amsfonts}
\usepackage{graphicx}
\usepackage{epstopdf}
\usepackage{algorithmic}
\Crefname{ALC@unique}{Line}{Lines}
\usepackage{graphicx}
\usepackage{epstopdf}
\usepackage{url}
\usepackage{amsmath}
\usepackage{mathtools}
\usepackage{subcaption}
\usepackage{booktabs}
\usepackage{multirow}
\usepackage{cite}
\usepackage{breqn}

\newsiamremark{example}{Example}
\newsiamthm{assumption}{Assumption}

\ifpdf
  \DeclareGraphicsExtensions{.eps,.pdf,.png,.jpg}
\else
  \DeclareGraphicsExtensions{.eps}
\fi

\numberwithin{theorem}{section}

\newcommand{\TheTitle}{A General Return-Mapping Framework for Fractional Visco-Elasto-Plasticity}
\newcommand{\TheAuthors}{Jorge Suzuki, Maryam Naghibolhosseini, Mohsen Zayernouri}

\headers{General Return-Mapping for Fractional Visco-Elasto-Plasticity}{\TheAuthors}

\title{{\TheTitle}}

\author{
  Jorge L. Suzuki\thanks{Department of Mechanical Engineering and Department
  of Computational Mathematics, Science and Engineering, Michigan State
  University, East Lansing, MI 48824, USA
  	(\email{suzukijo@msu.edu}).}
  \and
  Maryam Naghibolhosseini\thanks{Department of Communicative Sciences and Disorders, Michigan State
  University, East Lansing, MI 48824, USA
  	(\email{naghib@msu.edu}).}
  \and
  Mohsen Zayernouri\thanks{Department of Mechanical Engineering and Department
  of Statistics and Probability, Michigan State University, East Lansing, MI
  48824, USA
  	(\email{zayern@msu.edu}, \email{zayern@egr.msu.edu}), Corresponding Author.}
}

\usepackage{amsopn}

\ifpdf
\hypersetup{
  pdftitle={\TheTitle},
  pdfauthor={\TheAuthors}
}
\fi

\newsiamthm{remark}{Remark}

\DeclareMathOperator{\sign}{\mathrm{sign}}

\begin{document}

\maketitle

\begin{abstract}
We develop a fractional return-mapping framework for power-law visco-elasto-plasticity. In our approach, the fractional viscoelasticity is accounted through canonical combinations of Scott-Blair elements to construct a series of well-known fractional linear viscoelastic models, such as Kelvin-Voigt, Maxwell, Kelvin-Zener and Poynting-Thomson. We also consider a fractional quasi-linear version of Fung's model to account for stress/strain nonlinearity. The fractional viscoelastic models are combined with a fractional visco-plastic device, coupled with fractional viscoelastic models involving serial combinations of Scott-Blair elements. We then develop a general return-mapping procedure, which is fully implicit for linear viscoelastic models, and semi-implicit for the quasi-linear case. We find that, in the correction phase, the discrete stress projection and plastic slip have the same form for all the considered models, although with different property and time-step dependent projection terms. A series of numerical experiments is carried out with analytical and reference solutions to demonstrate the convergence and computational cost of the proposed framework, which is shown to be at least first-order accurate for general loading conditions. Our numerical results demonstrate that the developed framework is more flexible, preserves the numerical accuracy of existing approaches while being more computationally tractable in the visco-plastic range due to a reduction of $50\%$ in CPU time. Our formulation is especially suited for emerging applications of fractional calculus in bio-tissues that present the hallmark of multiple viscoelastic power-laws coupled with visco-plasticity.
\end{abstract}

\begin{keywords}
power-law visco-elasto-plasticity, time-fractional integration, fractional quasi-linear viscoelasticity.
\end{keywords}

\begin{AMS}
   34A08, 74A45, 74D10, 74S20, 74N30.
\end{AMS}

\section{Introduction}
Power law behavior has been observed in living cells \cite{Okajima2017, De2020} and bio-tissues \cite{Naghibolhosseini2018,Guo2021, suzuki2021bladder}. This stems from the ubiquitous self-similar and scale-free nature of the tissue/cell microstructure, which can be physically and mathematically scaled up to continuum level, manifesting in the power law behavior in the lumped sense. Such power law relationships have been seen in auditory hair cells, positioned in the sensory organ of hearing, cochlea \cite{Martin2001}, vocal fold tissues \cite{Choi2012b}, and bladder tissues \cite{suzuki2021bladder}.
Experimental evidence suggests that complex material behavior may possess more than a single power-law scaling in the viscoelastic regime, particularly in multi-fractal structures, such as in cells \cite{Stamenovic2007} and biological tissues \cite{vincent2012structural}, due to their complex, hierarchical and heterogeneous microstructures.For such cases, a single fractional rheological element is not sufficient to capture the observed behavior even if the data suggest a linear viscoelastic behavior. Stamenovi\'{c} et al. \cite{Stamenovic2007} measured the complex shear modulus of cultured human airway smooth muscle and observed two distinct power-law regimes, separated by an intermediate \textit{plateau}. Kapnistos \textit{et al.} \cite{Kapnistos2008} found an unexpected tempered power-law relaxation response of entangled polystyrene ring polymers, compared with the usual relaxation plateau of linear chain polymers. Such behavior was interpreted through self-similar conformations of double-folded loops in the ring polymers, instead of the reptation observed in linear chains.

In addition to multiple viscoelastic power-law behaviors, there also exists evidence of bio-plasticity in soft media \cite{Pajerowski2007,Bonadkar2016}. The creep behavior of human embryonic stem cells under differentiation was studied by Pajerowski \textit{et al.} \cite{Pajerowski2007} through micro-aspiration experiments at different pressures. The cell nucleus demonstrated distinguished visco-elasto-plastic power-law scalings, with $\alpha = 0.2$ for the plastic regime, independent of the applied pressure. It is discussed that such low power-law exponent arises due to the fractal arrangement of chromatin inside the cell nucleus. Studies on force-induced mechanical plasticity of mouse embryonic fibroblasts were performed by Bonadkar \textit{et al.} \cite{Bonadkar2016}. They found that the viscoelastic relaxation and the permanent deformations followed a stochastic, normally-distributed, power-law scaling $\beta(\omega)$ with values ranging from $\beta \approx 0$ to $\beta \approx 0.6$. The microstructural mechanism of plastic deformation in the cytoskeleton is due to the combination of permanent stretching and buckling of actin fibers.

Regarding existing modeling approaches of anomalous plasticity, several works employed fractional calculus to account for the visco-plastic regimes of different classes of materials \cite{suzuki2021fractional}. Three of the main approaches include: time-fractional, space-fractional and stress-fractional modeling. The time-fractional approaches focus on introducing memory effects into non-equilibrium viscous variables \cite{Suzuki2016,Xiao2017}, and consequently modeling power-laws in both viscoelastic and visco-plastic regimes, which is applied for polymers, cells, and tissues. Suzuki \textit{et al.} \cite{Suzuki2016} developed a fractional visco-elasto-plastic model that provides a constitutive interpolation between rate-independent plasticity and Perzyna's visco-plasticity by introducing a Scott-Blair (SB) model acting the plastic regime. This model utilizes a rate-dependent yield function, which was later proved to be thermodynamically consistent in a further extension of the model to account for continuum damage mechanics \cite{suzuki2021thermodynamically}. A three-dimensional space-fractional approach to elasto-plasticity was also developed by Sumelka \cite{sumelka2014application} in order to consider the spatial nonlocalities. The model is based on rate-independent elasto-plasticity, and nonlocal effects are modeled using a fractional continuum mechanics approach, where the strains are defined through a space-fractional Riesz-Caputo derivative of the displacements. Finally, the stress-fractional models for plasticity have found to be applicable for modeling the soil mechanics and geomaterials that follow a non-associated plastic flow \cite{Sumelka2014VP,Sumelka2019Soil}, \textit{i.e.} in which the yield surface expansion in the stress space does not follow the usual normality rule, and may be non-convex. Sumelka \cite{Sumelka2014VP} proposed a three-dimensional fractional visco-plastic model, where a fractional flow-rule with the order $0 < \alpha < 1$ in the stress domain naturally modeled the non-associative plasticity. This model recovers the classical Perzyna visco-plasticity as $\alpha \to 1$, and the effect of the fractional flow rule can be a compact descriptor of micro-structure anisotropy. Later on, Sun and Sumelka \cite{Sumelka2019Soil} developed a similar stress-fractional model, which was successfully applied for soils under compression. We refer the reader to the Sun \textit{et al.} review work on fractional calculus applications in plasticity \cite{sun2018new}.

In this work, we develop generalized fractional visco-elasto-plastic model, where the visco-plastic device can be coupled with several existing fractional linear/nonlinear viscoelastic representations. More specifically, we utilize a fractional visco-plastic device developed in \cite{Suzuki2016,suzuki2021thermodynamically}, which is then coupled with a series of linear fractional models, such as SB, Kelvin-Voigt (FKV), Maxwell (FM), Kelvin-Zener (FKZ), Poynting-Thomson (FPT); also a fractional quasi-linear viscoelastic (FQLV) model for large strains. Then, a generalized fractional return-mapping algorithm is proposed, which overcomes existing difficulties in previous developments by first fully discretizing all fractional operators, and then performing the predictor-corrector procedure. More specifically, the existing approaches are built on the notion of employing the predictor-corrector approach before the full discretization of fractional operators while treating trial states for stress and internal variables as continuous functions of time. This prevents models with serial combinations of SB elements to be incorporated in associated yield functions in a straightforward fashion. The main features of the proposed framework are:

\begin{itemize}
    \item We perform a full discretization of fractional viscoelastic models prior to the definition of trial states, which allows a linear decomposition between final and trial stresses regardless of the employed models.
    
    \item The fractional return-mapping algorithm is fully-implicit for linear viscoelastic rheology, and semi-implicit for quasi-linear viscoelasticity.
    
    \item Due to the full-discretization before the return-mapping procedure, the operations involving the plastic-slip are memoryless, which resembles return-mapping steps from the classical elasto-plasticity.
    
    \item The correction (return-mapping) step has the same structure regardless of the employed viscoelastic models.
\end{itemize}

We carry out a number of numerical experiments involving fabricated and reference solutions under monotone and general loading conditions, and observe a global accuracy ranging from $\mathcal{O}(\Delta t)$ to $\mathcal{O}(\Delta t^{2-\beta})$, according to the regularity induced by the associated fractional differential equations (FDEs) and loading conditions.

This work is organized as follows. In Section \ref{Sec:Definitions}, we present the mathematical definitions employed in this work. In Section \ref{Sec:FVE}, we describe the considered linear/quasi-linear fractional viscoelastic models, coupled with fractional visco-elasto-plasticity as explained in Section \ref{Sec:FVEP}. All corresponding models are discretized and posed in a unified fractional return-mapping form in Section \ref{Sec:RM}. Convergence analyses and computational performance evaluation of presented models and return-mapping algorithm are performed in Section \ref{Sec:NumericalResults}, followed by the concluding remarks in Section \ref{Sec:Conclusions}.

\section{Definitions of Fractional Calculus}
\label{Sec:Definitions}
%

We start with some preliminary definitions of fractional calculus
\cite{podlubny1998fractional}. The left-sided Riemann-Liouville integral of order {$\beta
	\in (0,1)$} is defined as
\begin{equation}
\label{Eq: left RL integral}
(\prescript{RL}{t_L}{\mathcal{I}}_{t}^{\beta} f) (t) = \frac{1}{\Gamma(\beta)}
\int_{t_L}^{t} \frac{f(s)}{(t - s)^{1-\beta} }\, ds,\,\,\,\,\,\, t>t_L,
\end{equation}
where $\Gamma$ represents the Euler gamma function and $t_L$ denotes the lower
integration limit. The corresponding inverse operator, the left-sided
fractional derivative of order $\beta$, is then defined based on (\ref{Eq: left
	RL integral}) as
\begin{equation}
\label{Eq: left RL derivative}
(\prescript{RL}{t_L}{\mathcal{D}}_{t}^{\beta} f) (t) = \frac{d}{dt}
(\prescript{RL}{t_L}{\mathcal{I}}_{t}^{1-\beta} f) (t) =
\frac{1}{\Gamma(1-\beta)} \frac{d}{dt} \int_{t_L}^{t} \frac{f(s)}{(t -
	s)^{\beta} }\, ds,\,\,\,\,\,\, t>t_L.
\end{equation}
{The left-sided Caputo derivative of order $\beta
	\in (0,1)$ is obtained as}
\begin{equation}
\label{Eq: left Caputo derivative}
(\prescript{C}{t_L}{\mathcal{D}}_{t}^{\beta} f) (t) =
(\prescript{RL}{t_L}{\mathcal{I}}_{t}^{1-\beta} \frac{df}{dt}) (t) =
\frac{1}{\Gamma(1-\beta)} \int_{t_L}^{t} \frac{f^{\prime}(s)}{(t - s)^{\beta}
}\, ds,\,\,\,\,\,\, t>t_L.
\end{equation}
The definitions of Riemann-Liouville and Caputo derivatives are linked by the following relationship:
	\begin{equation}
	\label{Eq: Caputo vs. Riemann}
	(\prescript{RL}{t_L}{\mathcal{D}}_{t}^{\beta} f) (t) =
	\frac{f(t_L)}{\Gamma(1-\beta) (t+t_L)^{\beta}} +
	(\prescript{C}{t_L}{\mathcal{D}}_{t}^{\beta} f) (t),
	\end{equation}
which can be obtained through integration by parts followed by the application of Leibniz rule on \eqref{Eq: left RL derivative}. We should note that the aforementioned derivatives coincide when dealing with homogeneous Dirichlet initial/boundary conditions. Finally, we define the two-parameter Mittag-Leffler function $E_{a,b}(z)$ \cite{Mainardi2011} as:
\begin{equation}
  E_{a,b}(z) = \sum^\infty_{k=0} \frac{z^k}{\Gamma(a k + b)},\quad Re(a) > 0,\quad b \in \mathbb{C},\quad z \in \mathbb{C}.
\end{equation}

\section{Fractional Viscoelasticity}
\label{Sec:FVE}

We present the linear and quasi-linear fractional viscoelastic models that we couple with the visco-plastic return-mapping procedure. 

\subsection{Linear Viscoelasticity}

\begin{figure}[!h]
     \centering
     \includegraphics[width=\columnwidth]{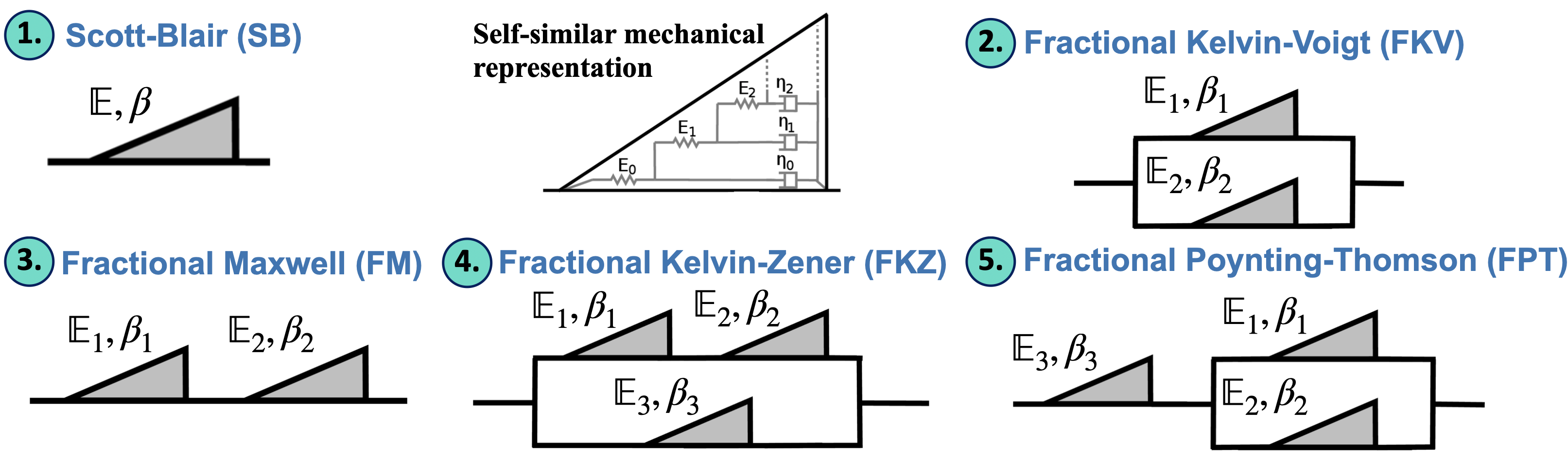}
     \caption{\label{fig:vediagram} Fractional linear viscoelastic models employed in this work, constructed from serial/parallel combinations of SB elements. The SB building blocks naturally account for an infinite fractal arrangement of Hookean/Newtonian elements. The employed fractional quasi-linear model is not represented by a mechanical analogue although the time-dependent component of the relaxation function has a SB-like representation.}
\end{figure}

\vspace{2mm}
\noindent\textbf{Scott-Blair (SB) Model:} The rheological \textit{building block} for our framework is the fractional SB viscoelastic element, which compactly represents an anomalous viscoelastic constitutive law relating the stresses and strains:
\begin{equation}\label{eq:SB}
\sigma(t) = \mathbb{E} \prescript{C}{0}{} \mathcal{D}_{t}^{\beta} \varepsilon(t), \quad t > 0, \quad \varepsilon(0) = 0,
\end{equation}
with pseudo-constant $\mathbb{E}_1\,[Pa.s^{\beta}] \ge 0$ and constant fractional order $0 < \beta < 1$, which provides a material interpolation between the Hookean ($\beta \to 0$) and Newtonian ($\beta \to 1$) elements. The pair $(\beta, \mathbb{E})$ uniquely represents the SB constants, where the \textit{pseudo-constant} $\mathbb{E}\,[Pa.s^{\beta}]$ compactly describes textural properties, such as the firmness of the material \cite{Blair1947, Jaishankar2013}. In this sense, $\mathbb{E}$ is interpreted as a snapshot of a non-equilibrium dynamic process instead of an equilibrium state. The corresponding rheological symbol for the SB model represents a fractal-like arrangement of springs and dashpots \cite{Schiessel1993,McKinley2013}, which we interpret as a compact, upscaled representation of a fractal-like microstructure. Regarding the thermodynamic admissibility, we refer the reader to Lion\cite{Lion1997} for the SB model, and Suzuki \textit{et al.}\cite{suzuki2021thermodynamically} for the combination of the SB element with more complex mechanisms of visco-plasticity and damage. The relaxation function $G(t)\,[Pa]$ for the SB model is given by the following inverse power-law form:
\begin{equation}\label{eq:SB_G}
G^{SB}(t) := \frac{\mathbb{E}}{\Gamma(1-\alpha)} t^{-\beta},
\end{equation}
which is the convolution kernel of the differ-integral form in (\ref{eq:SB}).

\vspace{2mm}
\noindent\textbf{Fractional Kelvin-Voigt (FKV) Model:} Through a parallel combination of SB elements, we obtain the following stress-strain relationship \cite{Schiessel1993}:
\begin{equation}\label{eq:FKV}
\sigma(t) = \mathbb{E}_1\, \prescript{C}{0}{} \mathcal{D}^{\beta_1}_t \varepsilon(t) + \mathbb{E}_2\, \prescript{C}{0}{} \mathcal{D}^{\beta_2}_t \varepsilon(t) ,\quad t > 0, \quad \varepsilon(0) = 0,
\end{equation}
with fractional orders $0 < \beta_1, \beta_2 < 1$ and associated pseudo-constants $\mathbb{E}_1\,[Pa.s^{\beta_1}] \ge 0$ and $\mathbb{E}_2\,[Pa.s^{\beta_2}] \ge 0$. The corresponding relaxation modulus $G(t)\,[Pa]$, is also an additive form of two SB elements:
\begin{equation}\label{eq:FKV_G}
G^{FKV}(t) := \frac{\mathbb{E}_1}{\Gamma(1-\beta_1)} t^{-\beta_1} + \frac{\mathbb{E}_2}{\Gamma(1-\beta_2)} t^{-\beta_2},
\end{equation}
which has a response characterized by two power-law regimes with a transition from faster to slower relaxation. Assuming $\beta_2 > \beta_1$, the asymptotic responses for small and large time-scales are given by $G^{FKV} \sim t^{-\beta_2}$ as $t\to 0$ and $G^{FKV} \sim t^{-\beta_1}$ as $t\to \infty$.

\vspace{2mm}
\noindent\textbf{Fractional Maxwell (FM) Model:} Through a serial combination of SB elements, we obtain the fractional Maxwell (FM) model\cite{Jaishankar2013}, given by:
\begin{equation}\label{eq:FM}
\sigma(t) + \frac{\mathbb{E}_2}{\mathbb{E}_1}\prescript{C}{0}{}\mathcal{D}^{\beta_2 - \beta_1}_t \sigma(t) = \mathbb{E}_2 \prescript{C}{0}{}\mathcal{D}^{\beta_2}_t \varepsilon(t) , \quad t>0,
\end{equation}
with pseudo-constants $\mathbb{E}_1\,[Pa.s^{\beta_1}] > 0$ and $\mathbb{E}_2 [Pa.s^{\beta_1}] \ge 0$, fractional orders $0 < \beta_1 < \beta_2 < 1$ with $0 < \beta_2 - \beta_1 < 1$ and two sets of initial conditions for strains $\varepsilon(0) = 0$ and stresses $\sigma(0) = 0$. We should note that in the case of non-homogeneous initial conditions, there needs to be a compatibility condition\cite{Mainardi2011} between stresses and strains at $t=0$. The corresponding relaxation function for this building block model assumes a more complex Miller-Ross form\cite{Jaishankar2013}:
\begin{equation}\label{eq:G_FMM}
G^{FM}(t) := \mathbb{E}_1 t^{-\beta_1} E_{\beta_2-\beta_1, 1-\beta_1}\left(-\frac{\mathbb{E}_1}{\mathbb{E}_2} t^{\beta_2 - \beta_1}\right).
\end{equation}
The presence of a Mittag-Leffler function in (\ref{eq:G_FMM}) leads to a stretched exponential relaxation for smaller times and a power-law behavior for longer times. We also observe that the limit cases are given by $G^{FM} \sim t^{-\beta_1}$ as $t\to 0$ and $G^{FM} \sim t^{-\beta_2}$ as $t \to \infty$, indicating that the FM model provides a behavior transitioning from slower-to-faster relaxation. We refer the reader to \cite{Jaishankar2013,bonfanti2020fractional,suzuki2021bladder} for a number of applications of the aforementioned models. We should notice that both FKV and FM models are able to recover the SB element with a convenient set of pseudo-constants and $\beta_1 = \beta_2$. 

\vspace{2mm}
\noindent\textbf{Fractional Kelvin-Zener (FKZ) model:} The fractional generalization of the standard linear solid (SLS) model is given by an FM branch in parallel with a third SB element, given by the following FDE:
\begin{equation}\label{eq:FKZ}
    \left[ 1 + \frac{\mathbb{E}_2}{\mathbb{E}_1} \prescript{C}{0}{} \mathcal{D}^{\beta_2-\beta_1}_t \right] \sigma(t) = \left[ \mathbb{E}_2\prescript{C}{0}{} \mathcal{D}^{\beta_2}_t + \mathbb{E}_3\prescript{C}{0}{} \mathcal{D}^{\beta_3}_t + \frac{\mathbb{E}_2 \mathbb{E}_3}{\mathbb{E}_1}\prescript{C}{0}{} \mathcal{D}^{\beta_2+\beta_3-\beta_1}_t \right] \varepsilon(t),
\end{equation}
with fractional orders $0<\beta_1 < \beta_2 < 1$ and conditions $0 < \beta_2 - \beta_1 < 1$ and $0 < \beta_2+\beta_3-\beta_1 < 1$, pseudo-constants $\mathbb{E}_1 \,[Pa.s^{\beta_1}] > 0$, $\mathbb{E}_2 \,[Pa.s^{\beta_2}] \ge 0$ and $\mathbb{E}_3 \,[Pa.s^{\beta_3}] \ge 0$, and the same initial conditions as in the FM model. We should note that the FM model is recovered when $\mathbb{E}_3 = 0$ and the FKV model is recovered when setting $\mathbb{E}_1 = 0$. The relaxation function is obtained in a straightforward fashion as the summation of the relaxation functions from the SB and FM models:
\begin{equation}\label{eq:FKZ_G}
    G^{FKZ}(t) := \mathbb{E}_1 t^{-\beta_1} E_{\beta_2-\beta_1, 1-\beta_1}\left(-\frac{\mathbb{E}_1}{\mathbb{E}_2} t^{\beta_2 - \beta_1}\right) + \frac{\mathbb{E}_3}{\Gamma(1-\beta_3)} t^{-\beta_3},
\end{equation}
which leads to three inverse power-law regimes for short, intermediate and long times with particular relationships between $\beta_1,\,\beta_2,\,\beta_3$ \cite{Schiessel1995}.

\vspace{2mm}
\noindent\textbf{Fractional Poynting-Thomson (FPT) Model:} Finally, we introduce our last fractional linear viscoelastic model, given by the serial combination between an FKV model and an SB element:
\begin{equation}\label{eq:FPT}
    \left[ 1 + \frac{\mathbb{E}_1}{\mathbb{E}_3} \prescript{C}{0}{} \mathcal{D}^{\beta_1-\beta_3}_t + \frac{\mathbb{E}_2}{\mathbb{E}_3} \prescript{C}{0}{} \mathcal{D}^{\beta_2-\beta_3}_t \right] \sigma(t) = \left[ \mathbb{E}_1\prescript{C}{0}{} \mathcal{D}^{\beta_1}_t + \mathbb{E}_2\prescript{C}{0}{} \mathcal{D}^{\beta_2}_t \right] \varepsilon(t),
\end{equation}
with $0 < \beta_3 < \beta_1 < 1$ and $0 < \beta_3 < \beta_2 < 1$, additional conditions $0 < \beta_1 - \beta_3 < 1$ and $0 < \beta_2 - \beta_3 < 1$, and pseudo-constants $\mathbb{E}_1\,[Pa.s^{\beta_1}] \ge 0$, $\mathbb{E}_2\,[Pa.s^{\beta_2}] \ge 0$ and $\mathbb{E}_3\,[Pa.s^{\beta_3}] > 0$ and homogeneous initial conditions $\sigma(0) = 0$ and $\varepsilon(0) = 0$. Similar to the FKZ model, we recover the FM model when setting either $\mathbb{E}_1$ or $\mathbb{E}_2$ to zero; although the FKV model cannot be recovered except for a trivial case when $\sigma(t) = 0$. 

\subsection{Quasi-Linear Fractional Viscoelasticity}

Although fractional linear viscoelastic models provide suitable relaxation functions that describe the anomalous viscoelastic dynamics of a number of soft materials, at times, complex microstructural deformation mechanisms and large strains induce material nonlinearities, hence, the relaxation function itself depends on the applied strain levels. To incorporate this additional effect, we also consider the following FQLV model \cite{fung2013biomechanics,Craiem2008}:
\begin{equation}\label{eq:QLV}
		\sigma(t,\varepsilon) = \int^t_0 G(t-s)\frac{\partial \sigma^e(\varepsilon)}{\partial \varepsilon} \dot{\varepsilon}\,ds,
\end{equation}
where the convolution kernel is given by a multiplicative decomposition of a reduced relaxation function $G(t)$ and an instantaneous nonlinear elastic tangent response with stress $\sigma^e$. In the work by Craiem \textit{et al.}\cite{Craiem2008}, the reduced relaxation function has a fractional Kelvin-Voigt-like form with one of the SB replaced with a Hookean element. Here, we assume a simpler rheology and adopt an SB-like reduced relaxation in the form:
\begin{equation}\label{eq:SBGreduced}
  G(t) = E t^{-\alpha} / \Gamma(1-\alpha)
\end{equation}
with the pseudo-constant $E$ with units $[s^\alpha]$. We adopt the same, two-parameter, exponential nonlinear elastic part as in\cite{Craiem2008}:
\begin{equation}\label{eq:expelastic}
  \sigma^{e}(\varepsilon) = A \left(e^{B \varepsilon} - 1\right),
\end{equation}
with $A$ having units of $[Pa]$. Plugging in \eqref{eq:SBGreduced} and \eqref{eq:expelastic} into \eqref{eq:QLV}, we obtain:
\begin{equation}\label{eq:FQLV}
		\sigma(t,\varepsilon) = \frac{E A B}{\Gamma(1-\alpha)}\int^t_0 \frac{e^{B \varepsilon(s)} \dot{\varepsilon}(s)}{(t-s)^\alpha}\,ds,
\end{equation}
which differs slightly from the linear SB model \eqref{eq:SB} in the sense that an additional exponential factor multiplies the function being convoluted.

\section{Fractional Visco-Elasto-Plasticity}
\label{Sec:FVEP}

\begin{figure}[!h]
     \centering
     \includegraphics[width=0.8\columnwidth]{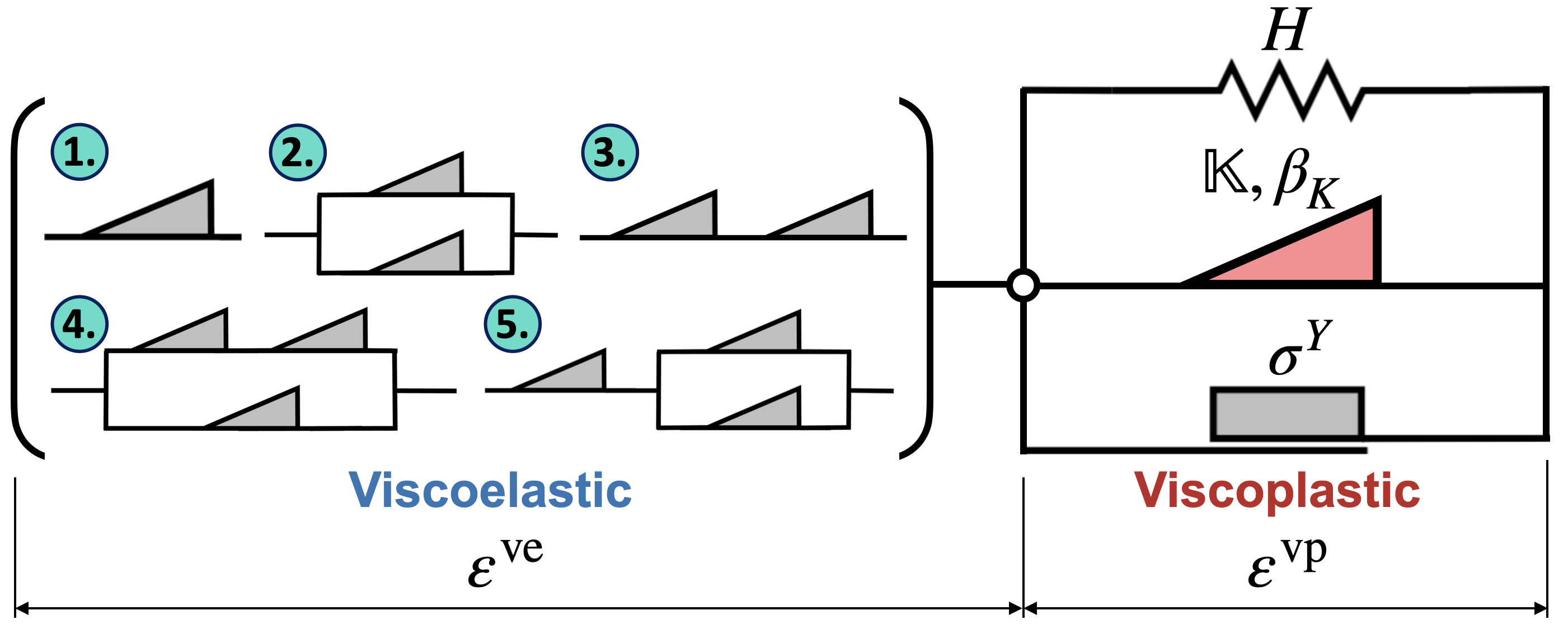}
     \caption{\label{fig:vpdiagram} Fractional visco-elasto-plastic model diagram. Here, any of the linear and quasi-linear fractional viscoelastic models can be separately coupled with a fractional visco-plastic rheological device.}
\end{figure}

With all fractional viscoelastic models defined in Section \ref{Sec:FVE}, we can couple any of them, subject to a viscoelastic strain $\varepsilon^{ve}(t)$, to the fractional visco-plastic device, illustrated in Fig. \ref{fig:vpdiagram}. The visco-plastic device is composed of a parallel combination of a Coulomb element with initial yield stress $\sigma^Y\,[Pa]$, a SB element with pseudo-constant $\mathbb{K}\,[Pa.s^{\beta_K}]$ and fractional order $\beta_K$, and a Hookean spring with constant $H\,[Pa]$. The entire visco-plastic part is subject to a visco-plastic strain $\varepsilon^{vp}(t):\mathbb{R}^+ \to \mathbb{R}$. In order to obtain the kinematic equations for the internal variables, we start with an additive decomposition of the total logarithmic strain $\varepsilon(t):\mathbb{R}^+ \to \mathbb{R}$ acting on the visco-elasto-plastic device:
\begin{equation}
    \varepsilon(t) = \varepsilon^{ve}(t) + \varepsilon^{vp}(t)
\end{equation}

The visco-plastic effects are accounted for through the definition of a memory- and rate-dependent yield function $f\left( \sigma, \alpha\right):\mathbb{R}\times\mathbb{R}^+\to\mathbb{R}^-\cup\lbrace
0\rbrace$ in the following form \cite{Suzuki2016}:
\begin{equation}\label{eq:Yield_function}
	f(\sigma, \alpha) := \lvert \sigma \rvert - \left[ \sigma^Y + \mathbb{K}\, 
	{}^C_0 \mathcal{D}^{\beta_K}_t \left(\alpha \right) + H\alpha \right].
\end{equation}
Here, $\alpha \in \mathbb{R}^+$ represents the internal hardening variable, and the above form accounts for the isotropic hardening. The set of admissible stresses lie in a closed convex space, where the associated boundary respects the yield condition of classical plasticity \textit{(see \cite{suzuki2021thermodynamically}, Lemma 4.1, setting the damage as $D = 0$)}. From the defined yield function (\ref{eq:Yield_function}) and the principle of maximum plastic dissipation \cite{Simo1998}, the following {properties hold:} \textbf{i)} associativity of the flow rule, 
\textbf{ii)} associativity in the hardening law, \textbf{iii)} Kuhn-Tucker complimentary conditions, and 
\textbf{iv)} convexity. The set of evolution equations for the internal variables $\varepsilon^{vp}$ and $\alpha$ is obtained by:
\begin{equation}
	\dot{\varepsilon}^{vp} = \frac{\partial f}{\partial \sigma} \dot{\gamma}, 
	\quad 
	\dot{\alpha} = -\frac{\partial f}{\partial R} \dot{\gamma},
\end{equation}
where $\dot{\gamma}(t):\mathbb{R}^+\to\mathbb{R}^+$ denotes the plastic slip rate. Evaluating the above equations using (\ref{eq:Yield_function}), we obtain the evolution for 
visco-plastic strains and hardening \cite{Suzuki2016}:
\begin{equation}\label{eq:evol_vp_damage}
\dot{\varepsilon}^{vp} = \mathrm{sign}(\tau) \dot{\gamma},
\end{equation}
\begin{equation}\label{eq:evol_hardening_damage}
\dot{\alpha} = \dot{\gamma}.
\end{equation}

\begin{proposition} \label{lem:closureslip}
The closure for the plastic slip rate $\dot{\gamma}(t) \in \mathbb{R}^+$ with a SB viscoelastic part of constants $(\mathbb{E},\,\beta_E)$, $(\mathbb{K},\,\beta_K)$, and $H$ (model M1 \cite{Suzuki2016}) with homogeneous initial conditions for the internal variables and their respective rates, \textit{i.e.}, $\varepsilon^{vp}(0) = \alpha(0) = \gamma(0) = 0$, $\dot{\gamma}(0) = 0$, and $\dot{\varepsilon}^{vp}(0) = \dot{\alpha}(0) = \dot{\gamma}(0) = 0$, is given by the following fractional Cauchy problem:
\begin{equation}
    \mathbb{E} \prescript{C}{0}{} \mathcal{D}^{\beta_E}_t \dot{\gamma}(t) 
    + \mathbb{K} \prescript{C}{0}{} \mathcal{D}^{\beta_K}_t \dot{\gamma}(t) + H \dot{\gamma}(t) = \sign(\sigma)\mathbb{E}\left[\frac{\dot{\varepsilon}(0) t^{-\beta_E}}{\Gamma(1-\beta_E)} + \prescript{C}{0}{} \mathcal{D}^{\beta_E}_t \dot{\varepsilon}(t) \right]
\end{equation}

\begin{proof}
		See Appendix \ref{Ap:dgamma_relaxation}.
\end{proof}
\end{proposition}


\section{A Class of Return-Mapping Algorithms for Fractional Visco-Elasto-Plasticity}
\label{Sec:RM}

Given the presented viscoelastic and visco-plastic models, respectively, in Sections \ref{Sec:FVE} and \ref{Sec:FVEP}, we now demonstrate how to solve each resulting system of nonlinear equations according to the choice of viscoelastic models. The considered fractional return-mapping approach in this work is fully discrete, \textit{i.e.}, we first discretize all fractional derivatives using a finite-difference approach, and then employ trial states for the internal variables in a predictor-corrector scheme.

We discretize the fractional Caputo derivatives in \eqref{eq:SB}-\eqref{eq:FM} through an implicit L1 finite-difference scheme\cite{lin2007finite}. Extensions to account for fast time-stepping approaches \cite{zeng2018stable} are straightforward, since they mostly affect the history terms computation. Let $\Omega = (0,T]$ be decomposed into a uniform time grid with $N$ time steps of size $\Delta t$, such that $t_n = n \Delta t$, with $n = 0,\,1,\,\dots,\,N$. The time-fractional Caputo derivative of a real-valued function $u(t) \in C^2(\Omega)$ at time $t = t_{n+1}$ is therefore discretized as \cite{lin2007finite}:
\begin{equation}\label{eq:L1}
    {}^C_0\mathcal{D}^\beta_t u(t)\big|_{t=t_{n+1}} = \frac{1}{\Delta t^\beta \Gamma(2-\beta)}\left[ u_{n+1} - u_n + \mathcal{H}^\alpha u \right] + \mathcal{O}(\Delta t^{2-\beta}),
\end{equation}
with the history term $\mathcal{H}^\nu u$ given by the following form:
\begin{equation}
\mathcal{H}^\beta u = \sum^n_{j=1} b^\beta_j \left[ u_{n+1-j} - u_{n-j}\right]
\end{equation}
with weights $b^\beta_j := (j+1)^{1-\beta} - j^{1-\beta}$. 

\subsection{Time-Fractional Integration of Viscoelastic Models}

In the following, we present the discretized forms for each considered fractional viscoelastic model from Section \ref{Sec:FVE}, which are represented in a fully-implicit fashion.
\vspace{2mm}

\noindent\textbf{Scott-Blair Model}: Evaluating both sides of \eqref{eq:SB} at $t = t_{n+1}$, we obtain:
\begin{equation}
    \sigma_{n+1} = \mathbb{E}_1\,{}^C_0\mathcal{D}^{\beta_1}_t \varepsilon(t)\big|_{t=t_{n+1}}
\end{equation}
in which by applying \eqref{eq:L1}, we directly obtain:
\begin{equation}\label{eq:SBdiscrete}
  \sigma_{n+1} = C^{SB}_1 \left[ \varepsilon_{n+1} - \varepsilon_{n} + \mathcal{H}^{\beta_1} \varepsilon \right]
\end{equation}
with the strain history $\mathcal{H}^{\beta_1} \varepsilon$ and constant $C^{SB}_1$, shown in Appendix \ref{Ap:constants} for the SB and the following model discretizations. 
\vspace{2mm}

\noindent\textbf{Fractional Kelvin-Voigt Model}: Evaluating both sides of \eqref{eq:FKV} at $t = t_{n+1}$, we obtain:
\begin{equation}
        \sigma_{n+1} = \mathbb{E}_1{}^C_0\mathcal{D}^{\beta_1}_t \varepsilon(t)\big|_{t=t_{n+1}} + \mathbb{E}_2{}^C_0\mathcal{D}^{\beta_2}_t \varepsilon(t)\big|_{t=t_{n+1}},
\end{equation}
which, applying \eqref{eq:L1}, for the fractional derivatives of order $\beta_1$ and $\beta_2$, leads to:
\begin{equation}\label{eq:FKVdiscrete}
  \sigma_{n+1} = C^{KV}_1 \left[ \varepsilon_{n+1} - \varepsilon_{n} + \mathcal{H}^{\beta_1} \varepsilon \right] + C^{KV}_2 \left[ \varepsilon_{n+1} - \varepsilon_{n} + \mathcal{H}^{\beta_2} \varepsilon \right].
\end{equation}

\noindent\textbf{Fractional Maxwell Model}: Evaluating both sides of \eqref{eq:FM} at $t = t_{n+1}$, we obtain:
\begin{equation}
    \sigma_{n+1} + \frac{\mathbb{E}_2}{\mathbb{E}_1}\prescript{C}{0}{}\mathcal{D}^{\beta_2 - \beta_1}_t \sigma(t)\big|_{t=t_{n+1}} = \mathbb{E}_2 \prescript{C}{0}{}\mathcal{D}^{\beta_2}_t \varepsilon(t)\big|_{t=t_{n+1}}
\end{equation}
in which applying \eqref{eq:L1} for the fractional derivatives of strains and stresses, leads to:
\begin{equation}\label{eq:FMdiscrete}
  \sigma_{n+1} = \frac{C^{M}_1 \left[ \varepsilon_{n+1} - \varepsilon_{n} + \mathcal{H}^{\beta_2} \varepsilon \right] + C^{M}_2 \left[ \sigma_{n} - \mathcal{H}^{\beta_2-\beta_1} \sigma \right]}{ 1 + C^{M}_2 },
\end{equation}
with the emergence of a stress history term $\mathcal{H}^{\beta_2-\beta_1} \sigma$.
\vspace{2mm}

\noindent\textbf{Fractional Kelvin-Zener Model}: Evaluating both sides of \eqref{eq:FKZ} at $t = t_{n+1}$, we obtain:
\begin{dmath*}
    \sigma_{n+1} + \frac{\mathbb{E}_2}{\mathbb{E}_1} \prescript{C}{0}{} \mathcal{D}^{\beta_2-\beta_1}_t\sigma(t)\big|_{t=t_{n+1}} =  \mathbb{E}_2\prescript{C}{0}{} \mathcal{D}^{\beta_2}_t \varepsilon(t)\big|_{t=t_{n+1}}+ \mathbb{E}_3\prescript{C}{0}{} \mathcal{D}^{\beta_3}_t \varepsilon(t)\big|_{t=t_{n+1}}+ \frac{\mathbb{E}_2 \mathbb{E}_3}{\mathbb{E}_1}\prescript{C}{0}{} \mathcal{D}^{\beta_2+\beta_3-\beta_1}_t \varepsilon(t)\big|_{t=t_{n+1}},
\end{dmath*}
which after applying \eqref{eq:L1} for the fractional derivatives of strains and stresses, leads to:
\begin{dmath}\label{eq:FKZdiscrete}
    \sigma_{n+1} = (1+C^{KZ}_4)^{-1} \left[ C^{KZ}_1 \left(  \Delta\varepsilon_{n+1}+\mathcal{H}^{\beta_2}\varepsilon\right)+C^{KZ}_2\left(\Delta\varepsilon_{n+1}+\mathcal{H}^{\beta_3}\varepsilon \right) + C^{KZ}_3\left(\Delta\varepsilon_{n+1}+\mathcal{H}^{\beta_2+\beta_3-\beta_1}\varepsilon\right)+C^{KZ}_4\left(\sigma_n-\mathcal{H}^{\beta_2-\beta_1}\sigma\right) \right]
\end{dmath}
with $\Delta \varepsilon_{n+1} = \varepsilon_{n+1} - \varepsilon_n$.
\vspace{2mm}

\noindent\textbf{Fractional Poynting-Thomson Model}: Finally, we evaluate both sides of \eqref{eq:FPT} and obtain:
\begin{dmath*}
    \sigma_{n+1} + \frac{\mathbb{E}_1}{\mathbb{E}_3} \prescript{C}{0}{} \mathcal{D}^{\beta_1-\beta_3}_t\sigma(t)\big|_{t=t_{n+1}} + \frac{\mathbb{E}_2}{\mathbb{E}_3} \prescript{C}{0}{} \mathcal{D}^{\beta_2-\beta_3}_t \sigma(t)\big|_{t=t_{n+1}} = \mathbb{E}_1\prescript{C}{0}{} \mathcal{D}^{\beta_1}_t \varepsilon(t)\big|_{t=t_{n+1}} + \mathbb{E}_2\prescript{C}{0}{} \mathcal{D}^{\beta_2}_t \varepsilon(t)\big|_{t=t_{n+1}},
\end{dmath*}
which after applying \eqref{eq:L1} for the fractional derivatives of strains and stresses, leads to:
\begin{dmath}\label{eq:FPTdiscrete}
    \sigma_{n+1} = (1 + C^{PT}_3 + C^{PT}_4)^{-1} \left[C_1\left(\Delta\varepsilon_{n + 1} + \mathcal{H}^{\beta_1}\varepsilon\right) + C^{PT}_2\left(\Delta\varepsilon_{n + 1} + \mathcal{H}^{\beta_2}\varepsilon \right) + C^{PT}_3\left(\sigma_{n + 1} + \mathcal{H}^{\beta_1 - \beta_3}\sigma\right) + C^{PT}_4\left(\sigma_n - \mathcal{H}^{\beta_2 - \beta_3}\sigma\right) \right].
\end{dmath}
\vspace{2mm}

\noindent\textbf{Fractional Quasi-Linear Viscoelastic Model:} The discretization for the FQLV model \eqref{eq:FQLV} has a slightly different development than the preceding models. It involves a slight modification of the fully-implicit L1 difference approach by a trapezoidal rule, taken on the exponential factor. More specifically, we evaluate the FQLV operator as:
\begin{equation}
    \sigma_{n+1} = \frac{E A B}{\Gamma(1-\beta)} \sum^n_{k=0} \int^{t_{k+1}}_{t_k} \left(t_{n+1}-s\right)^{-\beta} \exp(B \varepsilon_{k+\frac{1}{2}}) \left(\frac{\varepsilon_{k+1}-\varepsilon_k}{\Delta t}\right)\,ds
\end{equation}
with $\varepsilon_{i+\frac{1}{2}} = (\varepsilon_i + \varepsilon_{i+1})/2$. Following similar steps as in \cite{lin2007finite}, we obtain the following discretized stresses at $t = t_{n+1}$ for the FQLV model:
\begin{equation}\label{eq:FQLVdiscrete}
  \sigma_{n+1} = C^{QLV}_1\left[\exp(B \varepsilon_{n+\frac{1}{2}}) \left(\varepsilon_{n+1} - \varepsilon_n\right) + \mathcal{H}^\alpha \left(\varepsilon,\frac{\partial \sigma^e}{\partial\varepsilon}\right)\right]
\end{equation}
with constant $C^{QLV}_1 = E A B / (\Delta t^\alpha \Gamma(2-\alpha))$. The discretized history load in this case is given by:
\begin{equation}
  \mathcal{H}^\alpha \left(\varepsilon,\frac{\partial \sigma^e}{\partial\varepsilon}\right) = \sum^n_{k=1} \exp(B \varepsilon_{n-k+\frac{1}{2}}) \left(\varepsilon_{n-k+1}-\varepsilon_{n-k}\right) b_k
\end{equation}
with weights $b_k = (k+1)^{1-\alpha} - k^{1-\alpha}$. Since the trapezoid approximation of the strains in the exponential term are second-order accurate, the overall accuracy of the viscoelastic models is still bounded by the native L1-difference approach, and therefore, should be of $\mathcal{O}(\Delta t^{2-\alpha})$.

\begin{remark}
    We note that except for the FQLV model, any of the aforementioned discretizations for the linear models can recover the existing classical counterparts by properly setting $\beta_i \to 0$ or $\beta_i \to 1$. In these cases, to achieve a comparable performance to the integer-order models, history terms can be selectively disregarded and the corresponding discretization constants can be adjusted to their integer-order counterparts.
\end{remark}

\subsection{Time-Fractional Integration of Visco-Plasticity}

We start with the discretization of internal variables. Following \cite{Suzuki2016}, we assume a strain-driven process with known total strains $\varepsilon_{n+1}$ at time $t_{n+1}$. The strain decomposition becomes:
\begin{equation}\label{eq:discretizedkinematics}
    \varepsilon_{n+1} = \varepsilon^{ve}_{n+1} + \varepsilon^{vp}_{n+1}.
\end{equation}
The flow rule (\ref{eq:evol_vp_damage}) is discretized through a first-order backward-Euler approach, which yields:
\begin{equation}\label{eq:discretized_vp}
\varepsilon^{vp}_{n+1} = \varepsilon^{vp}_{n} + \sign(\sigma_{n+1})\Delta\gamma_{n+1}
\end{equation}
with $\Delta\gamma_{n+1} = \gamma_{n+1} - \gamma_n$, representing the plastic slip 
increment in the interval $[t_n,\,t_{n+1}]$. Similarly, the discretization of 
the hardening law (\ref{eq:evol_hardening_damage}) is given by
\begin{equation}\label{eq:discretized_alpha}
\alpha_{n+1} = \alpha_n + \Delta\gamma_{n+1}.
\end{equation}
Evaluating the yield function (\ref{eq:Yield_function}) at $t_{n+1}$ and employing discretization \eqref{eq:L1} for the hardening variable, we obtain:
\begin{align}\label{eq:discrete_f}
f_{n+1} & = \lvert\sigma_{n+1}\rvert - \left[\sigma^Y + 
\mathbb{K}\, {}^C_0 \mathcal{D}^{\beta_K}_t (\alpha)\big|_{t=t_{n+1}} + 
H\alpha_{n+1}\right] \nonumber \\
        & = \lvert\sigma_{n+1}\rvert - \left[\sigma^Y + 
\mathbb{K}^* \left(\alpha_{n+1} - \alpha_n + \mathcal{H}^{\beta_K}\alpha\right) +
H\alpha_{n+1}\right]
\end{align}
with $\mathbb{K}^* = \mathbb{K}/(\Delta t^{\beta_K}\Gamma(2-\beta_K))$.

The next step is to define the trial states for the stress and yield functions, which is the core idea to define the viscoelastic prediction phase, and the correction step after solving the internal visco-plastic variables. Therefore, we freeze the internal variables for the prediction step at $t_{n+1}$. Accordingly, the trial visco-plastic strains and hardening are given by:
\begin{equation}\label{eq:trialstates}
	\varepsilon^{vp^{trial}}_{n+1} = \varepsilon^{vp}_n, \quad 
	\alpha^{trial}_{n+1} = \alpha_n.
\end{equation}
In this token, the trial yield function is given by setting the above relationship for the hardening variable into \eqref{eq:discrete_f} to obtain:
\begin{equation}\label{eq:ftrial}
	f^{trial}_{n+1} = \lvert\sigma^{trial}_{n+1}\rvert - \left[\sigma^Y + \mathbb{K}^* \left(\mathcal{H}^{\beta_K}\alpha\right) + H\alpha_{n} \right].
\end{equation}

In order to complete the return-mapping procedure, we need an explicit relationship between the stresses $\sigma_{n+1}$ in terms of the known total strains $\varepsilon_{n+1}$. To achieve this, we solve for the plastic slip $\Delta \gamma$ using a discrete consistency condition $f_{n+1} = 0$. We start with the trial stresses for each presented fractional viscoelastic model by substituting the visco-plastic trial strain \eqref{eq:trialstates} and \eqref{eq:discretizedkinematics} into \eqref{eq:SBdiscrete}-\eqref{eq:FQLVdiscrete}, where we obtain the following for each discretized model:
\vspace{2mm}

\noindent\emph{Scott-Blair:}
\begin{equation}
  \sigma^{trial}_{n+1} = C^{SB}_1 \left[ \varepsilon_{n+1} - \varepsilon_{n} + \mathcal{H}^{\beta_1} (\varepsilon-\varepsilon^{vp}) \right]
\end{equation}

\noindent\emph{Fractional Kelvin-Voigt:}
\begin{equation}
  \sigma^{trial}_{n+1} = C^{KV}_1 \left[ \varepsilon_{n+1} - \varepsilon_{n} + \mathcal{H}^{\beta_1} (\varepsilon - \varepsilon^{vp}) \right] + C^{KV}_2 \left[ \varepsilon_{n+1} - \varepsilon_{n} + \mathcal{H}^{\beta_2} (\varepsilon - \varepsilon^{vp}) \right]
\end{equation}

\noindent\emph{Fractional Maxwell:}
\begin{equation}
  \sigma^{trial}_{n+1} = \frac{C^{M}_1 \left[ \varepsilon_{n+1} - \varepsilon_{n} + \mathcal{H}^{\beta_2} (\varepsilon - \varepsilon^{vp}) \right] + C^{M}_2 \left[ \sigma_{n} - \mathcal{H}^{\beta_2-\beta_1} \sigma \right]}{ 1 + C^{M}_2 }
\end{equation}

\noindent\emph{Fractional Kelvin-Zener:}
\begin{dmath}
    \sigma^{trial}_{n+1} = (1+C^{KZ}_4)^{-1} \left[ C^{KZ}_1 \left( \Delta\varepsilon_{n+1}+\mathcal{H}^{\beta_2}(\varepsilon-\varepsilon^{vp})\right)+C^{KZ}_2\left(\Delta\varepsilon_{n+1}+\mathcal{H}^{\beta_3}(\varepsilon-\varepsilon^{vp}) \right) + C^{KZ}_3\left(\Delta\varepsilon_{n+1}+\mathcal{H}^{\beta_2+\beta_3-\beta_1}(\varepsilon-\varepsilon^{vp})\right)+C^{KZ}_4\left(\sigma_n-\mathcal{H}^{\beta_2-\beta_1}\sigma\right) \right]
\end{dmath}

\noindent\emph{Fractional Poynting-Thomson:}
\begin{dmath}
    \sigma^{trial}_{n+1} = (1 + C^{PT}_3 + C^{PT}_4)^{-1} \left[C^{PT}_1\left(\Delta\varepsilon_{n + 1} + \mathcal{H}^{\beta_1}(\varepsilon-\varepsilon^{vp})\right) + C^{PT}_2\left(\Delta\varepsilon_{n + 1} + \mathcal{H}^{\beta_2}(\varepsilon-\varepsilon^{vp}) \right) + C^{PT}_3\left(\sigma_{n + 1} + \mathcal{H}^{\beta_1 - \beta_3}\sigma\right) + C^{PT}_4\left(\sigma_n - \mathcal{H}^{\beta_2 - \beta_3}\sigma\right) \right]
\end{dmath}

\noindent\emph{Fractional Quasi-Linear Viscoelastic Model:}

For this model, we follow a similar procedure of substituting the viscoelastic strains into \eqref{eq:FQLVdiscrete}, however, we evaluate the exponential term explicitly in time for all stages of the return-mapping algorithm. Therefore, the corresponding trial state becomes:
\begin{equation}
  \sigma^{trial}_{n+1} = C^{QLV}_1\left[\exp(B (\varepsilon_n - \varepsilon^{vp}_n)) \left(\varepsilon_{n+1} - \varepsilon_n\right) + \mathcal{H}^\alpha \left(\varepsilon-\varepsilon^{vp},\frac{\partial \sigma^e}{\partial\varepsilon}\right)\right].
\end{equation}

\subsection{Generalized Fractional Return-Mapping Algorithm}

From the aforementioned trial states, each discretized viscoelastic constitutive laws \eqref{eq:SBdiscrete}-\eqref{eq:FQLVdiscrete} and recalling \eqref{eq:discretized_vp}, one can show the following stress correction onto the yield surface:
\begin{equation}\label{eq:stressprojgeneral}
    \sigma_{n+1} = \sigma^{trial}_{n+1} - \sign(\sigma^{trial}_{n+1}) C^{ve}_{RM}(\mathbb{E},\,\Delta t,\,\varepsilon) \Delta\gamma_{n+1},
\end{equation}
where all discretized aforementioned viscoelastic models change the return-mapping procedure by a scaling factor $C^{ve}_{RM}(C,\,\varepsilon_n,\,\varepsilon^{vp}_n) \in \mathbb{R}^+$ acting on the Lagrange multiplier $\Delta \gamma$, which is given by the following for each model:
\begin{equation}\label{eq:projconstants}
    C^{ve}_{RM} = 
    \begin{cases}
        C^{SB}_1 & \emph{(Scott-Blair)} \\[5pt]
        C^{KV}_1 + C^{KV}_2 & \emph{(Fractional Kelvin-Voigt)} \\[5pt]
        C^{M}_1/(1 + C^{M}_2) & \emph{(Fractional Maxwell)} \\[5pt]
        (C^{KZ}_1 \!+\! C^{KZ}_2 \!+\! C^{KZ}_3) / (1\!+\!C^{KZ}_4) & \emph{(Fractional Kelvin-Zener)} \\[5pt]
        (C^{PT}_1 + C^{PT}_2)/(1+C^{PT}_3+C^{PT}_4) & \emph{(Fractional Poynting-Thomson)} \\[5pt]
        C^{QLV}_1 \exp(B(\varepsilon_n - \varepsilon^{vp}_n)) & \emph{(Fractional Quasi-Linear-Viscoelastic)}
    \end{cases}
\end{equation}

We show the derivation of \eqref{eq:stressprojgeneral} for the fractional Kelvin-Zener model in Appendix \ref{Ap:RMFKZ}, from which the Scott-Blair, fractional Maxwell and fractional Kelvin-Voigt models can be directly recovered; note that the derivation for the fractional Poynting-Thomson and quasi-linear viscoelasticity follow similarly in a straightforward fashion. Substituting the updated stresses \eqref{eq:stressprojgeneral} into the discrete yield function \eqref{eq:discrete_f} and recalling \eqref{eq:ftrial}, we obtain:
\begin{equation}
    f_{n+1} = f^{trial}_{n+1} - \left( C^{ve}_{RM} + \mathbb{K}^* + H\right) \Delta \gamma.
\end{equation}
Enforcing the discrete yield condition $f_{n+1} = 0$, we obtain the solution for the discrete plastic slip:
\begin{equation}\label{eq:discreteslip}
\boxed{
    \Delta \gamma_{n+1} = \frac{f^{trial}_{n+1}}{C^{ve}_{RM} + \mathbb{K}^* + H}
}.
\end{equation}

\begin{algorithm}[!ht]
	\caption{Fractional return-mapping algorithm.}
	\label{alg:RMFMM}
	\begin{algorithmic}[1]
		\STATE{Database for $\varepsilon$, $\varepsilon^{vp}$, $\sigma$, $\alpha$, and total strain $\varepsilon_{n+1}$.}
		\STATE{$\varepsilon^{vp^{trial}}_{n+1} = \varepsilon^{vp}_n, \quad
			\alpha^{trial}_{n+1} = \alpha_n$}
		\STATE{Compute $\sigma^{trial}_{n+1}$ from \eqref{eq:SBdiscrete}-\eqref{eq:FQLVdiscrete} according to the selected fractional viscoelastic model.}
		\STATE{$f^{trial}_{n+1} = \lvert\sigma^{trial}_{n+1}\rvert - \left[\sigma^Y + \mathbb{K}^* \left(\mathcal{H}^{\beta_K}\alpha\right) + H\alpha_{n} \right]$}
		\IF{$f^{trial}_{n+1} \le 0$}
		\STATE{$\varepsilon^{vp}_{n+1} = \varepsilon^{vp}_n$,\quad
			$\alpha_{n+1} = \alpha_n$,\quad $\sigma_{n+1} =
			\sigma^{trial}_{n+1}$.}
		\ELSE
		\STATE{Return-Mapping:}
		\STATE{Compute $C^{ve}_{RM}$ from \eqref{eq:projconstants} according to the selected fractional viscoelastic model.}
		\STATE{{$\Delta\gamma_{n+1} = f^{trial}_{n+1}/(C^{ve}_{RM} + K^* + H)$}}
		\STATE{{$\sigma_{n+1} = \sigma^{trial}_{n+1} - \sign(\sigma^{trial}_{n+1}) C^{ve}_{RM} \Delta\gamma$}}
		\STATE{$\varepsilon^{vp}_{n+1} = \varepsilon^{vp}_{n} +
			\sign(\tau_{n+1})\Delta\gamma$}
		\STATE{$\alpha_{n+1} = \alpha_n + \Delta\gamma$}
		\ENDIF

	\end{algorithmic}
\end{algorithm}

\subsubsection{Comparison of the Return-Mapping Algorithm to the Existing Approaches}

In \cite{Suzuki2016}, trial states were defined prior to the discretization of fractional operators, and the corresponding trial variables were taken as continuous functions of time, therefore, making the return-mapping procedure ``semi-discrete." Let the quantities with bars, $(\bar{\cdot})$, be the corresponding solutions for the procedure developed in \cite{Suzuki2016}. For the SB viscoelastic case, one has the following trial stresses at $t = t_{n+1}$:
\begin{equation}
    \bar{\sigma}^{\text{trial}}_{n+1} = \mathbb{E} \prescript{C}{0}{}\mathcal{D}^{\beta_E}_t (\varepsilon - \bar{\varepsilon}^{{vp}^{trial}})\big|_{t=t_{n+1}}
\end{equation}
in which, after employing the discretized plastic flow rule, the following relationship between the corrected and trial stresses is obtained:
\begin{equation}
    \bar{\sigma}_{n+1} = \bar{\sigma}^{\text{trial}}_{n+1} - \mathbb{E}\sign(\bar{\sigma}_{n+1}) \prescript{C}{0}{}\mathcal{D}^{\beta_E}_t (\Delta \gamma) \big|_{t=t_{n+1}}.
\end{equation}
This equation can be explicitly be inserted into the discrete yield function to solve for the plastic slip rate. While such procedure is straightforward for SB and FKV viscoelastic elements, it is non-trivial for the serial combinations such as the FM, FKZ and FPT models. For instance, if we follow the same procedure for the FM model, we obtain:
\begin{equation}
\bar{\sigma}^{trial}_{n+1} + \frac{\mathbb{E}_2}{\mathbb{E}_1}\prescript{C}{0}{}\mathcal{D}^{\beta_2 - \beta_1}_t (\bar{\sigma}^{trial})\big|_{t=t_{n+1}} = \mathbb{E}_2 \prescript{C}{0}{}\mathcal{D}^{\beta_2}_t (\varepsilon-\bar{\varepsilon}^{{vp}^{trial}})\big|_{t=t_{n+1}},
\end{equation}
which yields the following relationship between $\bar{\sigma}$ and $\bar{\sigma}^{trial}$:
\begin{equation}
\left( \bar{\sigma}_{n+1} - \bar{\sigma}^{trial}_{n+1}\right) + \frac{\mathbb{E}_2}{\mathbb{E}_1}\prescript{C}{0}{}\mathcal{D}^{\beta_2 - \beta_1}_t (\bar{\sigma} - \bar{\sigma}^{trial})\big|_{t=t_{n+1}} = -\mathbb{E}_2 \prescript{C}{0}{}\mathcal{D}^{\beta_2}_t (\varepsilon-\bar{\varepsilon}^{{vp}^{trial}})\big|_{t=t_{n+1}}.
\end{equation}
Except for the SB case, a fractional viscoelastic model involving a serial combination of SB elements cannot be incorporated into the yield function in a differential form, unless a full discretization is performed at this stage. This happens since the discretized yield function \eqref{eq:discrete_f} requires a closed description of $\sigma_{n+1}$, needing an equivalent Boltzmann representation for such models, which is impractical due to complex forms of relaxation kernels. Therefore, our approach in this work already carries the trial states with fully-discretized fractional operators, which closely and completely resembles classical elastoplastic approaches.

Regarding the obtained discretizations in this work, we note that the plastic slip \eqref{eq:discreteslip} assumes a simple form similar to the rate-independent elasto-plasticity. As discussed above, in the return-mapping procedure, developed in \cite{Suzuki2016}, the trial states and plastic slip were assumed to have memory in the discretization procedure, and therefore, a fractional relaxation equation in the following form was obtained:
\begin{equation}\label{eq:slip_M1}
    \Delta \bar{\gamma}_{n+1} = \frac{\mathbb{E}^*\left(\Delta \bar{\gamma}_n - \mathcal{H}^{\beta_E}\Delta \bar{\gamma} \right) + \mathbb{K}^*\left(\Delta \bar{\gamma}_n - \mathcal{H}^{\beta_K}\Delta \bar{\gamma} \right) + \bar{f}^{trial}_{n+1}}{\mathbb{E}^* + \mathbb{K}^* + H}.
    \end{equation}
Furthermore, we observe that the obtained plastic slip discretization in this work has two less history terms to be evaluated. Although this does not influence the computational complexity of the original scheme, we show in the numerical examples that this fact still leads to about $50\,\%$ less CPU time. Regarding the difference in the stress solutions, let $t = t_p$ be the time step of onset of plasticity for the first time. Therefore, we have the following estimate:
\begin{equation}\label{eq:diff_slip}
    \lvert \sigma_{p+1} - \bar{\sigma}_{p+1} \rvert = \frac{\mathbb{E}^*}{\mathbb{E}^* + \mathbb{K}^* + H} \left[\mathbb{K}^*\left(\mathcal{H}^{\beta_E}\Delta\bar{\gamma} - \mathcal{H}^{\beta_K}\Delta\bar{\gamma} \right) - H\left(\Delta \bar{\gamma}_p - \mathcal{H}^{\beta_E} \Delta \bar{\gamma}\right) \right],
\end{equation}
which shows that at such stage, both discretizations coincide when $\beta_E = \beta_K$ and $H = 0$. In the following Section, we verify such estimate by obtaining an analytical solution with the aid of Proposition \ref{lem:closureslip}.

\section{Numerical Tests}
\label{Sec:NumericalResults}

We present three convergence examples with different loading conditions to verify the employed fractional viscoelastic models, the validity of the new fractional visco-plastic return-mapping algorithm, and the full visco-elasto-plastic response of the models. For the convergence 
	analyses, let $\mathbf{u}^*$ and $\mathbf{u}^\delta$ be, respectively, the 
	reference and approximate solutions in $\Omega =(0,T]$, for a specific 
	time-step size $\Delta t$. We define the following relative error measures:
\begin{equation}\label{eq:errordef}
\mathrm{err}_N(\Delta t) = \frac{\vert u^*_N \!-\! u^\delta_N 
	\vert}{\vert u^*_N \vert}, \,\,
	\mathrm{err}(\Delta t) = \frac{\vert \vert \mathbf{u}^* \!-\! \mathbf{u}^\delta 
	\vert \vert_{L^2(\Omega)}}{\vert \vert \mathbf{u}^* \vert 
	\vert_{L^2(\Omega)}}, \,\, \mathrm{Order} = \log_2 
\left[\frac{\textrm{err}(\Delta t)}{\mathrm{err}(\Delta t / 2)}\right].
\end{equation}

We consider homogeneous initial conditions for all model variables in all 
cases. The presented algorithms were implemented in MATLAB R2020b and were run 
in a system with Intel Core i7-8850H CPU with 2.60 GHz, 32 GB RAM and MacOS 11.5 operating system.

\begin{example}[Convergence of fractional viscoelastic algorithms]

We perform a convergence study of the fractional viscoelastic component of our framework under the stress relaxation and monotone loading experiments. For this example, we set $(\mathbb{E}_1,\mathbb{E}_2,\mathbb{E}_3) = (1,\,1,\,1)$ and $(\beta_1,\,\beta_2,\,\beta_3) = (0.3,\,0.7,\,0.1)$ for fractional linear viscoelastic models, ensuring all fractional derivatives are taken with an equivalent order $\beta \in (0,1)$. For the FQLV model, we set $E = 1$, $\beta = 0.3$, $A = 1$, and $B = 1$.

For the stress relaxation test, we impose a step strain $\varepsilon(t) = H(t) \varepsilon_0$ with $\varepsilon_0 = 1$ for $T = 1000\,[s]$, where $H(t)$ denotes the Heaviside step function. We compare the obtained solutions at $t = T$ for the SB, FKV, FM, and FKZ models to their corresponding relaxation functions \eqref{eq:SB_G}, \eqref{eq:FKV_G}, \eqref{eq:G_FMM}, and \eqref{eq:FKZ_G}. The FPT and QLV models are not analyzed in this step since their time-dependent stress relaxation functions are not readily available, and they are instead analyzed under the monotone strains. The obtained results are illustrated in Fig. \ref{fig:ex1_rel_error}, where an expected linear convergence behavior is obtained for all models, given the non-smooth nature of the stress relaxation solution.

\begin{figure}[!h]
    \centering
     \begin{subfigure}[b]{0.49\columnwidth}
         \centering
         \includegraphics[width=\textwidth]{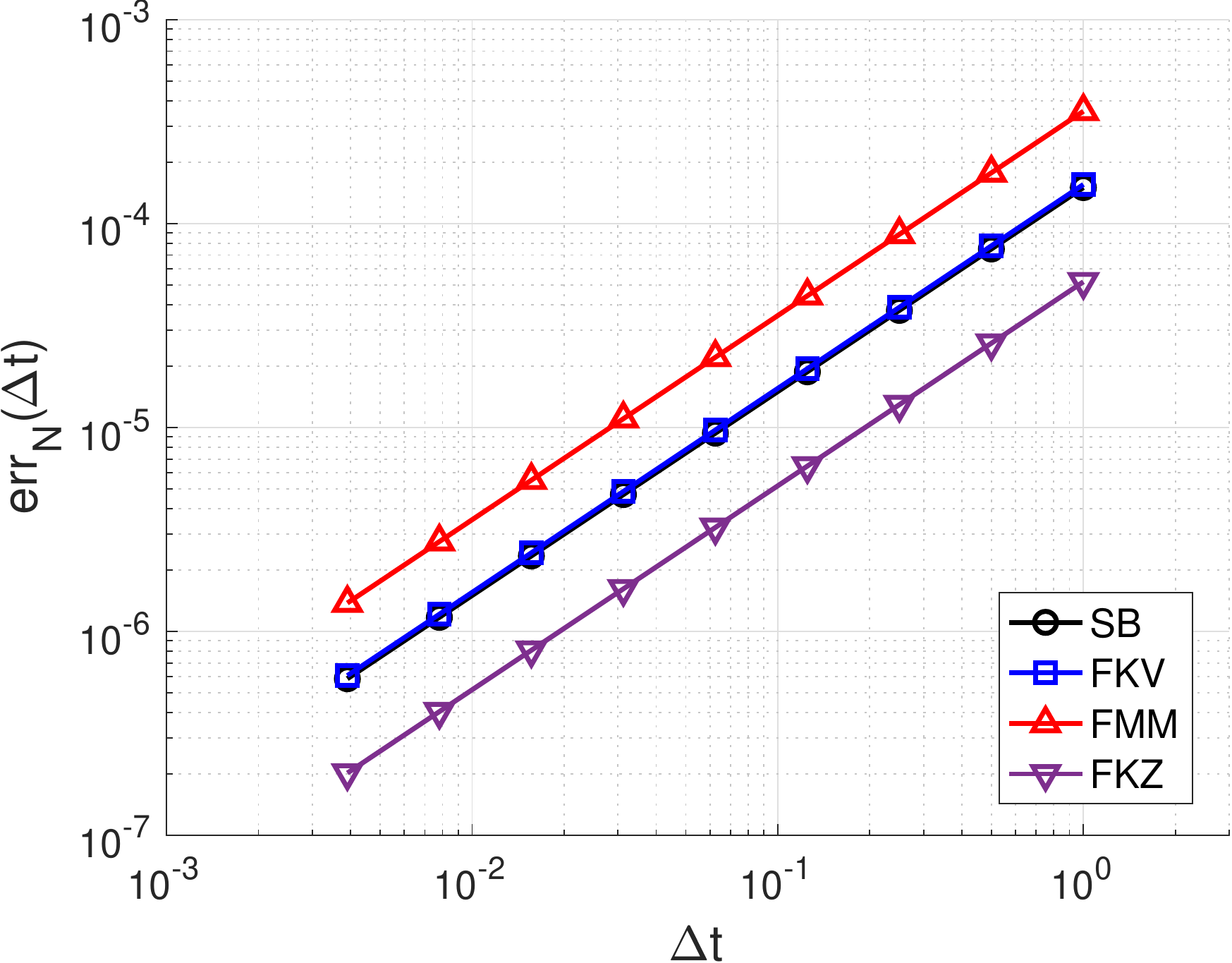}
         \caption{Stress relaxation test.}
         \label{fig:ex1_rel_error}
     \end{subfigure}
     \hfill
     \begin{subfigure}[b]{0.49\columnwidth}
         \centering
         \includegraphics[width=\textwidth]{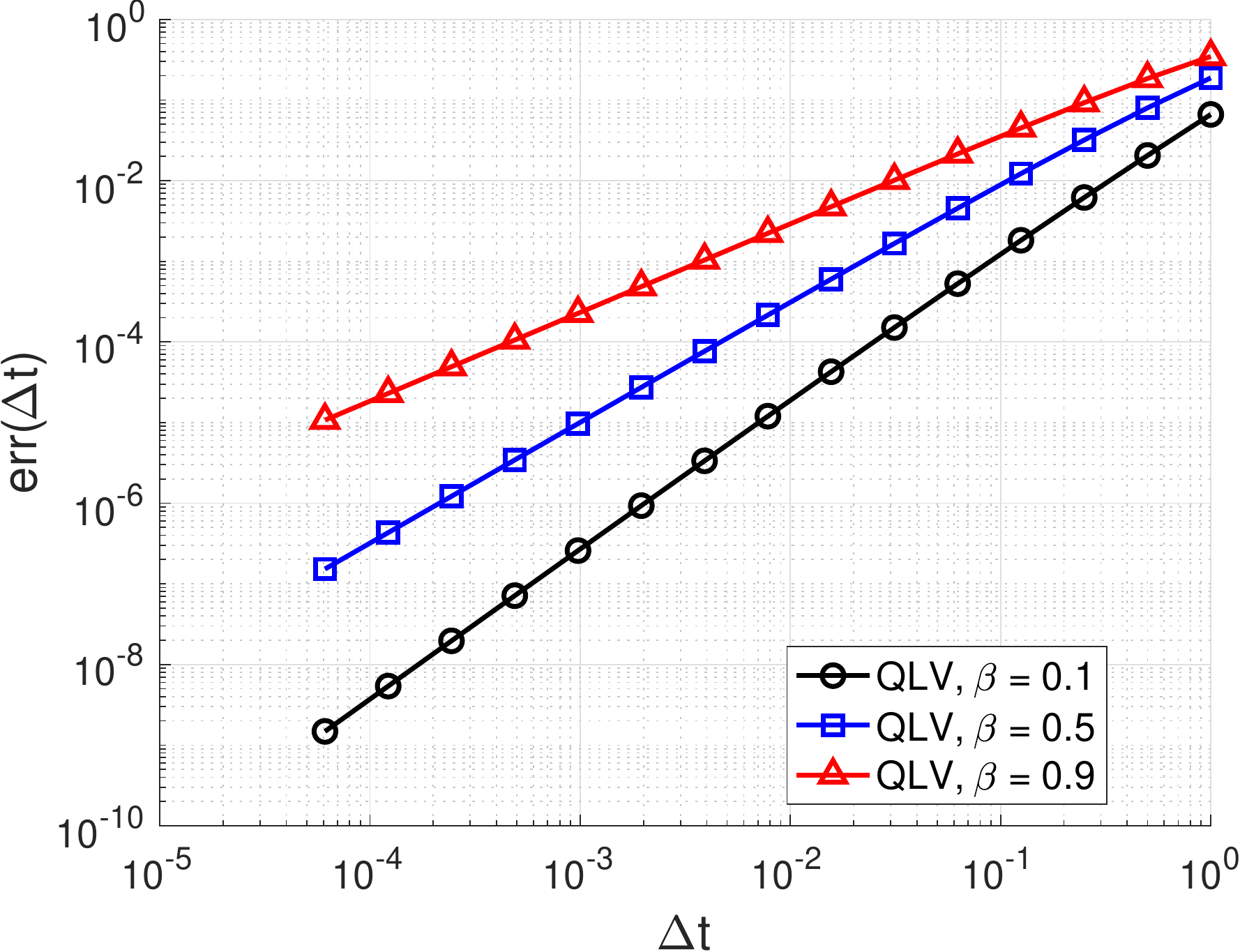}
         \caption{Monotone strain test - FQLV.}
         \label{fig:ex1_monotone_error_QLV}
     \end{subfigure}
     \caption{\label{fig:ex1_conv} Convergence analysis for the fractional viscoelastic models with known analytical solutions. \textit{(a)} A stress relaxation test with non-smooth step-strains, and material parameters $(\mathbb{E}_1, \mathbb{E}_2, \mathbb{E}_3) = (1,1,1)$ and $(\beta_1, \beta_2, \beta_3) = (0.3, 0.7, 0.1)$, yielding first-order convergence. \textit{(b)} Convergence for the FQLV model with a fabricated solution of linearly increasing strains and material properties $(E,\,\beta,\,A,\,B) = (1, 0.3, 1, 1)$. The slopes of the error curves are $q \approx 2-\beta$.}
\end{figure}

For the monotone strain case, we set $T=1$ and fabricate a solution for strains in the form $\varepsilon(t) = \varepsilon_T (t/T)$ with the total applied strain fixed at $\varepsilon_f = 1$. Since analytical solutions for all fractional viscoelastic models are difficult to obtain, we compute a reference solution for each model with $\Delta t = 2^{-17}\,[s]$. Particularly for the FQLV model, we utilize the the fabricated strain function $\varepsilon(t)$ to obtain the following analytical stress solution:
\begin{equation}
    \sigma^{{QLV}^*}(t) = E A B^\beta \exp(B t)\left[ 1 - \frac{\Gamma(1-\beta, B t)}{\Gamma(1-\beta)}\right],
\end{equation}
where $\Gamma(\cdot,\cdot)$ denotes the upper incomplete gamma function. The convergence results for all fractional viscoelastic models with respect to the reference numerical solution are presented in Fig. \ref{fig:ex1_monotone_ref}, while the results for the FQLV model with the analytical solution are illustrated in Fig. \ref{fig:ex1_monotone_error_QLV}. We observe for both cases that the accuracy of the implemented and developed schemes is of order $\mathcal{O}(\Delta t^{2-\beta})$. The difference in error slopes among models in Fig. \ref{fig:ex1_monotone_ref} is due to the highest fractional order assigned to each model. For the SB and FQLV models, the fractional order is set as $\beta = 0.3$, and therefore, the observed slope is $q \approx 1.7$. For all remaining models and choice of fractional orders, the error slopes are determined by the fractional derivative of the highest order, which is $\beta_2 = 0.7$ in this example, yielding $q \approx 1.3$.

\begin{figure}[!h]
     \centering
     \includegraphics[width=0.49\columnwidth]{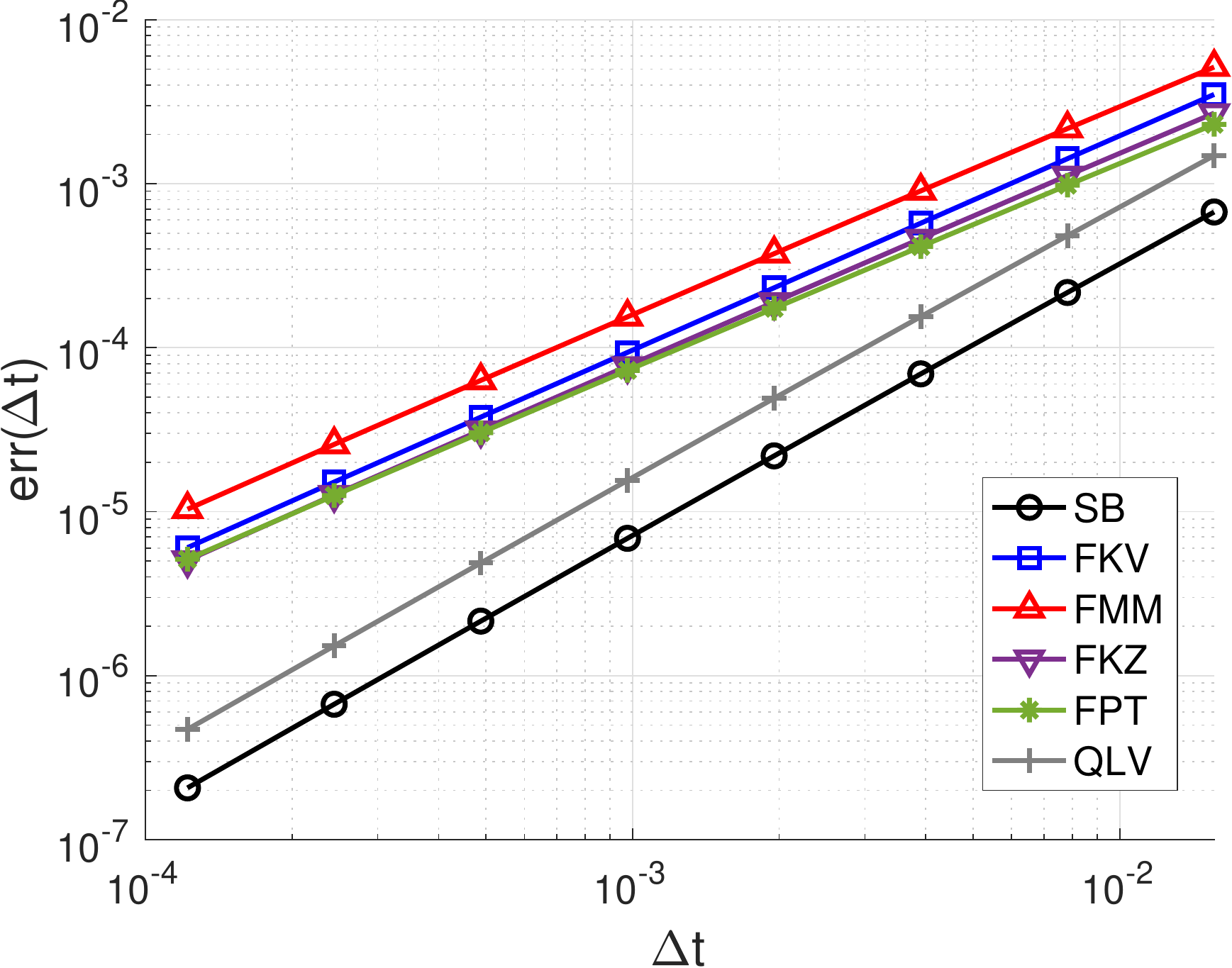}
     \caption{\label{fig:ex1_monotone_ref} Convergence analysis for all fractional viscoelastic models with $(\mathbb{E}_1, \mathbb{E}_2, \mathbb{E}_3) = (1,1,1)$ and $(\beta_1, \beta_2, \beta_3) = (0.3, 0.7, 0.1)$. A cubic strain function was employed with a reference solution using the time step size $\Delta t = 2^{-17}$. Monotone loading test with the convergence rate of $q \approx 1.3$ was applied for all models.}
\end{figure}

\end{example}

\begin{example}[Convergence of fractional visco-plastic algorithms]

The purpose of this example is to demonstrate the conditions where the presented plastic slip discretization \eqref{eq:discreteslip}, the form \eqref{eq:slip_M1} from \cite{Suzuki2016} and their associated return-mapping algorithms are equivalent, and also provide a numerical estimate for their difference when such conditions are not satisfied. For this purpose, we test a monotone load where an analytical solution is available, and a case with a cyclic load under high strain rates. For both cases, we set a SB viscoelastic part with $\mathbb{E} = 50\,[Pa.s^{\beta_E}]$ and $\mathbb{K} = 5\,[Pa.s^{\beta_K}]$.

For the monotone strain case, we start with a fabricated solution for strains in the form $\varepsilon(t) = A t^3$ with $A = \varepsilon_T/T^3\,[s^{-3}]$. Here, $\varepsilon_f$ denotes the total applied stress and $T$ represents the final simulation time. Utilizing the result of Lemma \ref{lem:closureslip} and setting $\beta_E = \beta_K = 0$ and $\sigma^Y = H = 0$, we obtain the following analytical solution for stresses:
\begin{equation}\label{eq:ex2_fabricated}
    \sigma^*(t) = \frac{6\, A\, \mathbb{E}\, \mathbb{K}}{\mathbb{E}+\mathbb{K}} \frac{t^{3-\beta_E}}{\Gamma(4-\beta_E)}.
\end{equation}
We should note that the proposed fabricated solution ensures no internal variable is a linear function, and therefore, not computed exactly by the L1 discretization. We set $\varepsilon_T = 1$ and $T = 1\,[s]$, and therefore, $A = 1\,[s^{-3}]$. Table \ref{tab:Conv_Ex2} presents the obtained convergence results for the fabricated solution \eqref{eq:ex2_fabricated} for both return-mapping algorithms and under the same fractional-orders $\beta_E$ and $\beta_K$. We observe that the errors coincide for this particular case while the accuracy of order $\mathcal{O}(\Delta t^{2-\beta})$ of the L1 approach is also achieved. The computational times are illustrated in Fig.\ref{fig:ex1_time}, where the developed fractional return-mapping approach, when using a SB viscoelastic element, is about $50\%$ faster than the original return-mapping approach from \cite{Suzuki2016} since about half of the history terms need to be computed.

\begin{table}[!htbp]
	\centering
	\caption{\label{tab:Conv_Ex2}Convergence behavior for the return-mapping Algorithm \ref{alg:RMFMM}, obtained in this work, and the original approach from \cite{Suzuki2016} for an FVEP device with an SB element.}
	\begin{tabular}{@{}lllllllll@{}}
		\toprule
		{}&\multicolumn{2}{c}{$\beta_E=\beta_K = 0.1$}&{}&\multicolumn{2}{c}{$\beta_E=\beta_K = 0.5$}&{}&\multicolumn{2}{c}{$\beta_E=\beta_K = 0.9$} \\
		\cline{2-3}\cline{5-6}\cline{8-9}
		$\Delta t$ & $\mbox{err}(\Delta t)$ & $\mbox{Order}$ & {} & $\mbox{err}(\Delta t)$ & $\mbox{Order}$ & {} & $\mbox{err}(\Delta t)$ & $\mbox{Order}$ \\ \midrule
		$2^{-9}$ & 3.2426e--06 & -- & {} & 9.2971e--05 & -- & {} & 1.3246e--03 & -- \\
        $2^{-10}$ & 9.1853e--07 & 1.8197 & {} & 3.3109e--05 & 1.4895 & {} & 6.1875e--04 & 1.0981 \\
        $2^{-11}$ & 2.5845e--07 & 1.8294 & {} & 1.1763e--05 & 1.4929 & {} & 2.8884e--04 & 1.0991 \\
        $2^{-12}$ & 7.2323e--08 & 1.8374 & {} & 4.1731e--06 & 1.4951 & {} & 1.3479e--04 & 1.0995 \\
        $2^{-13}$ & 2.0145e--08 & 1.8440 & {} & 1.4788e--06 & 1.4966 & {} & 6.2895e--05 & 1.0998 \\
        $2^{-14}$ & 5.5891e--09 & 1.8497 & {} & 5.2369e--07 & 1.4977 & {} & 2.9344e--05 & 1.0999 \\
        \bottomrule
    \end{tabular}
\end{table}

\begin{figure}[!ht]
    \centering
    \includegraphics[width=0.49\columnwidth]{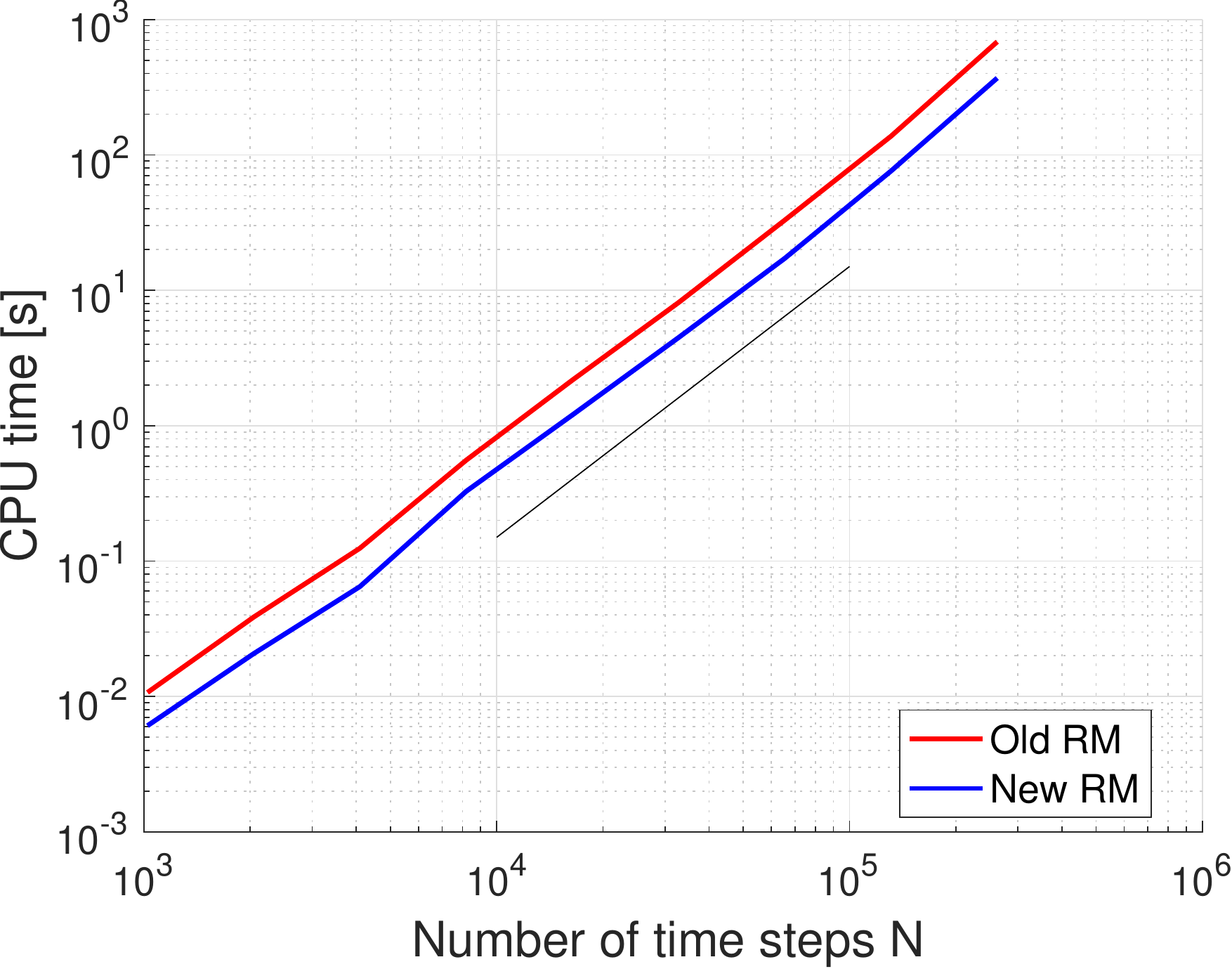}
    \caption{\label{fig:ex1_time}CPU times for the developed fractional return-mapping algorithm and the original one \cite{Suzuki2016} for an SB viscoelastic part. The black line has slope $q=2$.}
\end{figure}

Similar results are obtained for the monotone loading condition, however, this is not the case under general loadings. To demonstrate the difference between the visco-elasto-plastic discretization $\sigma_{n+1}$, developed in this work, and $\bar{\sigma}_{n+1}$ from \cite{Suzuki2016}, we take the latter as a reference solution with $\Delta t = 2^{-19}\,[s]$ and $T=1\,[s]$. We also consider $\sigma^Y = 10\,[Pa]$, $\beta_E = 0.3$ and $\beta_K = 0.7$ with the same pseudo-constants as in the previous test case. A constant rate loading/unloading cyclic strain test of the following form is employed:
\begin{equation}\label{eq:cyclic_strain}
	\varepsilon(t) = \frac{2 \varepsilon_A}{\pi} 
	\arcsin\left(\sin\left(2\pi\omega 
	t\right)\right),
\end{equation}
where we consider a strain amplitude $\varepsilon_A = 0.25$, and two strain frequencies of $\omega = 1\,[Hz]$ and $\omega = 60\,[Hz]$. The difference between both approaches is illustrated in Fig.\ref{fig:ex2_hysteresis}. Here, higher frequencies result in higher strain rates, and consequently a significant plastic strain history even after a number of hysteresis cycles. The obtained results confirm the estimates from \eqref{eq:diff_slip}, which is already valid at the onset of plasticity. Furthermore, we observe that a tenfold increase in strain rates approximately leads to a tenfold increase in the difference between both algorithms.

\begin{figure}[!ht]
    \centering
     \begin{subfigure}[b]{0.49\columnwidth}
         \centering
         \includegraphics[width=\textwidth]{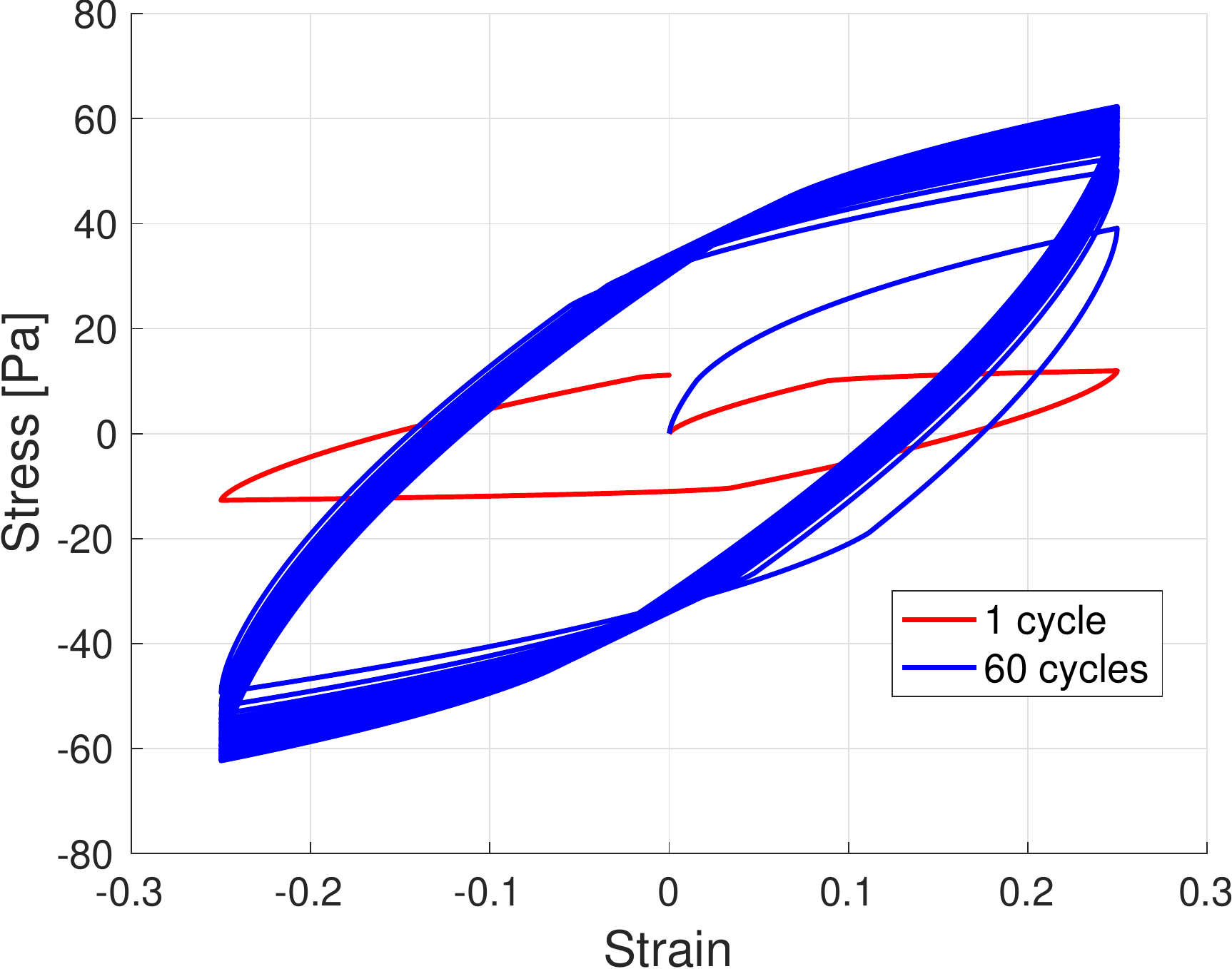}
         \caption{Reference solution.}
         \label{fig:1cycle}
     \end{subfigure}
     \hfill
     \begin{subfigure}[b]{0.49\columnwidth}
         \centering
         \includegraphics[width=\textwidth]{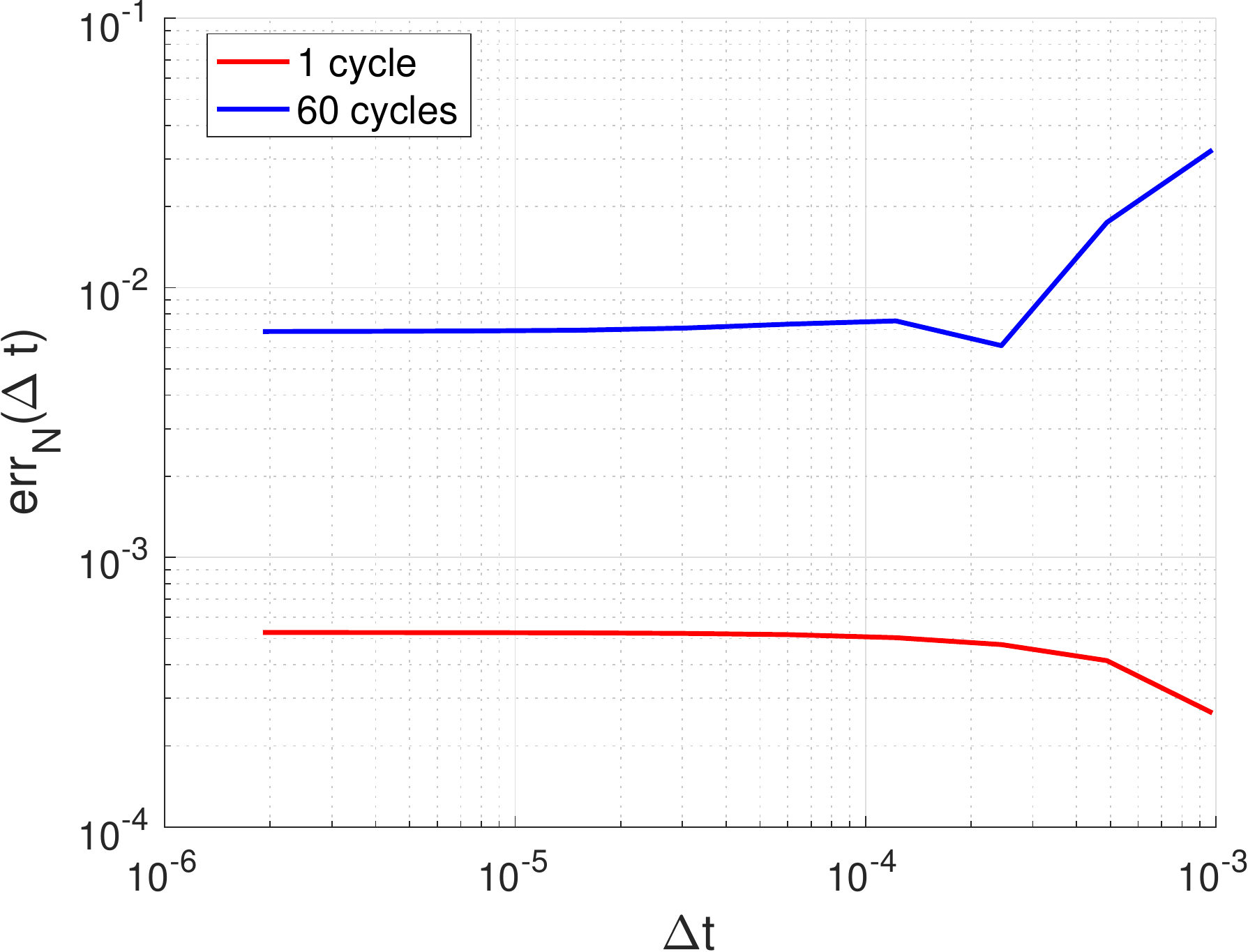}
         \caption{Endpoint error between both approaches.}
         \label{fig:60cycles}
     \end{subfigure}
     \caption{\label{fig:ex2_hysteresis} Comparison between the presented return-mapping algorithm and the reference approach from \cite{Suzuki2016} under low and high frequency loading.}
\end{figure}

\end{example}

\begin{example}[Convergence of fractional visco-elasto-plasticity]

Finally, we perform a verification on the entire fractional visco-elasto-plastic framework under cyclic strain. Since no fabricated solutions are available, we employ reference solutions with time step size $\Delta t = 2^{-18}\,[s]$. Let $T=1\,[s]$ with the same applied strains \eqref{eq:cyclic_strain} as in the previous example. The viscoelastic material parameters are set to $(\mathbb{E}_1, \mathbb{E}_2, \mathbb{E}_3) = (50,50,50)$ and $(\beta_1, \beta_2, \beta_3) = (0.3, 0.7, 0.1)$, in addition, the visco-plastic parameters are taken as $\mathbb{K}=5$, $\beta_K = 0.7$, $H=0$ and $\sigma^Y = 1$. Fig. \ref{fig:ex3} illustrates the obtained convergence results, where all models except the FKV one showed a convergence rate of order $q \approx 1.3$, which is compatible with the employed L1 discretization scheme and given $\beta_2 = \beta_K = 0.7$. The FKV model achieved linear asymptotic convergence for the considered example, which is the expected worst case scenario from the backward-Euler discretization of internal variables. We believe the difference in convergence behavior between the FKV model and the others could be due to the sharper response of the FKV model because of the stiffer rheology combined with the nonlinear loading/unloading response. This combination of effects could result in a lower solution regularity and therefore, a lower convergence rate.

\begin{figure}[!ht]
     \centering
     \includegraphics[width=0.49\columnwidth]{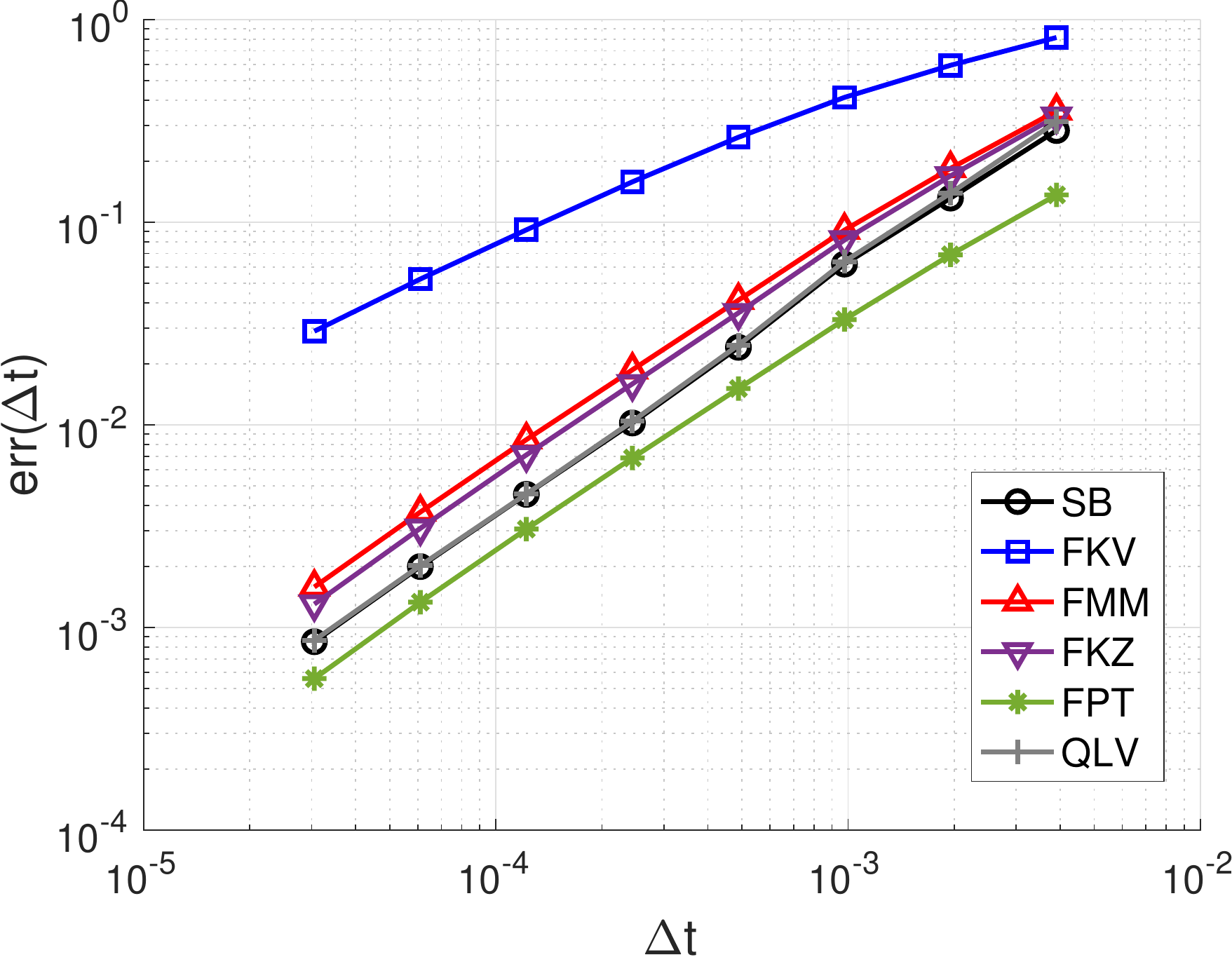}
     \caption{\label{fig:ex3} Convergence analysis for the fractional visco-elasto-plastic models under cyclic loads. Due to the particular choice of fractional orders (with $\beta_2 = \beta_K = 0.7$ being dominant), we observed the convergence rate off $q \approx 1.3$ for all models except for the FKV. In the latter case, we observe a linear convergence to the reference solution.}
\end{figure}

\begin{figure}[!ht]
    \centering
     \begin{subfigure}[b]{0.49\columnwidth}
         \centering
         \includegraphics[width=\textwidth]{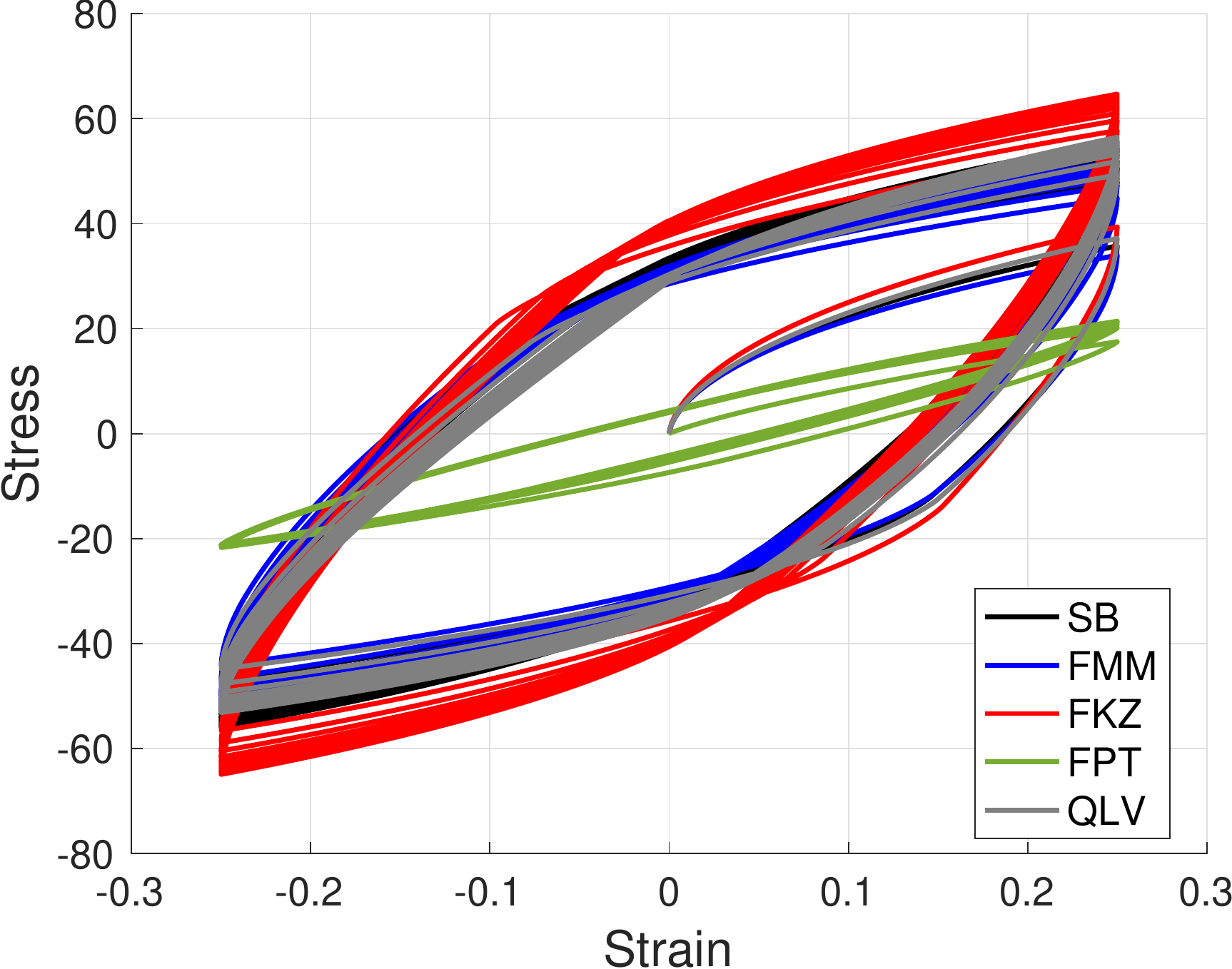}
         \caption{VEP reference solution.}
     \end{subfigure}
     \hfill
     \begin{subfigure}[b]{0.49\columnwidth}
         \centering
         \includegraphics[width=\textwidth]{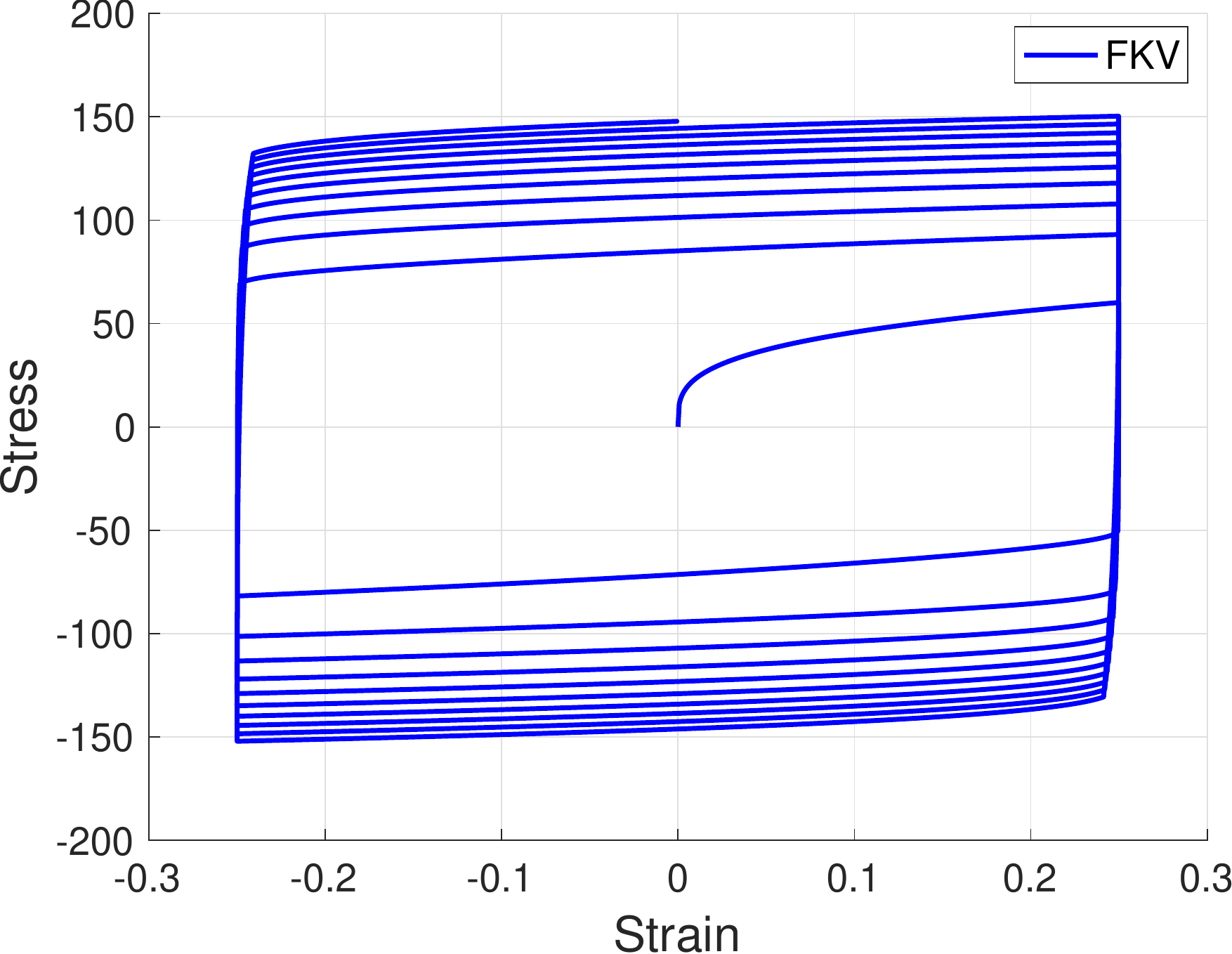}
         \caption{VEP reference solution - FKV model.}
     \end{subfigure}
     \caption{\label{fig:ex3_hysteresis} Visco-elasto-plastic reference solutions for the employed models for the first 30 loading cycles. We notice a similar behavior for most models under the choice of material parameters except for the FPT and FKV models. The FKV particularly yielded a very stiff response due to the combination of high fractional order values and high strain rates.}
\end{figure}

\end{example}

\section{Conclusions}
\label{Sec:Conclusions}

We proposed a general return-mapping procedure for multiple power-law, time-fractional visco-elasto-plastic materials. The developed framework provided a flexible way to integrate multiple known fractional viscoelastic models that are representative of soft materials rheology to power-law visco-plastic hardening and permanent strains. Furthermore, a nonlinear viscoelasticity, suitable for bio-tissues, was considered through a fractional quasi-linear Fung's model, which allowed the possibility of plasticity onset after substantial amounts of viscoelastic strains. The main features of the proposed framework are:

\begin{itemize}
    \item The trial states were taken after full discretization of stress and internal variables, which allowed a straightforward decomposition of the yield function into the final and trial states.
    
    \item The developed return-mapping procedure is fully implicit for linear viscoelastic models and semi-implicit for quasi-linear viscoelasticity. For simplicity, the chosen numerical discretization for fractional derivatives was an L1 finite-difference approach.
    
    \item Our correction step for visco-plasticity had the same structure for all viscoelastic models with the only difference being a scaling discretization constant.
    
    \item We carried out numerical experiments with analytical and reference solutions that demonstrated the $\mathcal{O}(\Delta t^{2-\beta})$ global accuracy, surprisingly even in some instances with general loading/unloading conditions.
    
    \item The developed return-mapping discretization was compared to an existing approach, and the difference between discretizations relied on cases with extensive plastic history and high strain rates.
\end{itemize}

Regarding the computational costs, the framework is computationally tractable since it does not involve history calculations for the plastic slip, and is therefore about $50\%$ faster than the existing fractional frameworks, regardless of the employed numerical discretization for fractional derivatives. Extensions on fast numerical schemes of order $\mathcal{O}(N \log N)$ for the employed time-fractional derivatives would be straightforward to implement. We also note that the thermodynamics of all models in the developed framework can be analyzed through the approach developed in \cite{suzuki2021thermodynamically}.

The modeling framework developed here could be applied to simulate the self-similar structures and memory-dependent behavior in human bio-tissues \cite{Naghibolhosseini2015, Naghibolhosseini2018, Suzuki2021}. The visco-elasto-plastic characteristics can be observed in different bio-tissues in human body specifically due the process of aging. Aging results in the oxidation or loss of elastin, which leads to the loss of tissue elasticity such as in the vocal fold tissues \cite{Branco2015, Ferster2017}. Furthermore, in terms of multi-scale modeling, the lumped plastic behavior introduced here could potentially be coupled with existing discrete dislocation dynamics (DDD) models \cite{de2021atomistic,de2022nonlocal}. The models developed in this work uniquely qualify for simulating such characteristics, which will be done in our future work.

\appendix

\section{Proof of Proposition \ref{lem:closureslip}} \label{Ap:dgamma_relaxation}

\begin{proof}
Similar to the derivation of the tangent elasto-plastic modulus in classical plasticity \cite{Simo1998}, we start by taking the time derivative of the yield function to enforce the persistency condition:
\begin{align}
    \dot{f}(\sigma,\alpha) & = \frac{d}{dt}\Bigl\{ \lvert \sigma(t) \rvert - \left[ \sigma^Y + \mathbb{K}\, {}^C_0 \mathcal{D}^{\beta_K}_t \alpha(t) + H\alpha(t) \right] \Bigr\} \\
    & = \sign(\sigma)\dot{\sigma}(t) - \left[\mathbb{K} \frac{d}{dt}{}^C_0 \mathcal{D}^{\beta_K}_t \alpha(t) + H\dot{\alpha}(t)\right]. \nonumber
\end{align}
Using the SB stress-strain relationship \eqref{eq:SB}, we obtain:
\begin{equation}\label{eq:dyield}
    \dot{f}(\sigma,\alpha) = \sign(\sigma)\mathbb{E} \left[\frac{d}{dt} {}^C_0 \mathcal{D}^{\beta_E}_t \varepsilon(t) \!-\! \frac{d}{dt} {}^C_0 \mathcal{D}^{\beta_E}_t \varepsilon^{vp}(t)\right] \!-\! \left[\mathbb{K} \frac{d}{dt}{}^C_0 \mathcal{D}^{\beta_K}_t \alpha(t) \!+\! H\dot{\alpha}(t)\right].
\end{equation}
Employing definition \eqref{Eq: left Caputo derivative} for the Caputo derivative, performing integration by parts and employing the Leibniz integral rule, we obtain:
\begin{align}\label{eq:dcaputo}
    \frac{d}{dt} {}^C_0 \mathcal{D}^{\beta}_t u(t) & = \frac{1}{\Gamma(1-\beta)} \frac{d}{dt} \int^t_0 \frac{\dot{u}(s)}{(t-s)^{\beta}}ds \quad \mathrm{(from \eqref{Eq: left Caputo derivative})} \nonumber \\
    & = \frac{1}{\Gamma(1-\beta)} \frac{d}{dt} \left[\dot{u}(s)\frac{(t-s)^{1-\beta}}{1-\beta}\big|^0_t + \int^t_0 \frac{(t-s)^{1-\beta}\ddot{u}(s)}{1-\beta}ds \right] \nonumber \\
    & = \frac{\dot{u}(0)t^{-\beta}}{\Gamma(1-\beta)} + \frac{1}{\Gamma(1-\beta)}\int^t_0 \frac{\ddot{u}(s)}{(t-s)^\beta}ds \nonumber \\
    & = \frac{\dot{u}(0)t^{-\beta}}{\Gamma(1-\beta)} + {}^C_0 \mathcal{D}^{\beta}_t \dot{u}(t).
\end{align}
Substituting \eqref{eq:dcaputo} into \eqref{eq:dyield}, setting $\dot{\gamma}(0) = 0$, and therefore $\dot{\alpha}(0) = 0$ from \eqref{eq:evol_hardening_damage} and $\dot{\varepsilon}^{vp}(0) = 0$ from \eqref{eq:evol_vp_damage}, we obtain:
\begin{equation}\label{eq:dyield1}
    \dot{f}(\sigma,\alpha) = \sign(\sigma)\mathbb{E} \left[ \frac{\dot{\varepsilon}(0) t^{-\beta_E}}{\Gamma(1-\beta_E)} + {}^C_0 \mathcal{D}^{\beta_E}_t \dot{\varepsilon}(t) - {}^C_0 \mathcal{D}^{\beta_E}_t \dot{\varepsilon}^{vp}(t)\right] - \mathbb{K}{}^C_0 \mathcal{D}^{\beta_K}_t \dot{\alpha}(t) - H\dot{\alpha}(t). 
\end{equation}
Finally, substituting \eqref{eq:evol_hardening_damage} and \eqref{eq:evol_vp_damage} into \eqref{eq:dyield1}, and enforcing the persistency condition $\dot{f}(\sigma,\alpha) = 0$, we obtain:
\begin{equation}
    \mathbb{E} \prescript{C}{0}{} \mathcal{D}^{\beta_E}_t \dot{\gamma}(t) 
    + \mathbb{K} \prescript{C}{0}{} \mathcal{D}^{\beta_K}_t \dot{\gamma}(t) + H \dot{\gamma}(t) = \sign(\sigma)\mathbb{E}\left[\frac{\dot{\varepsilon}(0) t^{-\beta_E}}{\Gamma(1-\beta_E)} + \prescript{C}{0}{} \mathcal{D}^{\beta_E}_t \dot{\varepsilon}(t) \right].
\end{equation}

\end{proof}

\section{Discretization Constants and Terms for Fractional Viscoelastic Models} \label{Ap:constants}
\vspace{2mm}

\noindent\emph{Scott-Blair:}
\begin{equation}
    C^{SB}_1 = \frac{\mathbb{E}}{\Delta t^{\beta_1} \Gamma(2-\beta_1)}
\end{equation}

\noindent\emph{Fractional Kelvin-Voigt:}
\begin{equation}
    C^{KV}_1 = \frac{\mathbb{E}_1}{\Delta t^{\beta_1} \Gamma(2-\beta_1)},\quad C^{KV}_2 = \frac{\mathbb{E}_2}{\Delta t^{\beta_2} \Gamma(2-\beta_2)}
\end{equation}

\noindent\emph{Fractional Maxwell:}
\begin{equation}
    C^M_1 = \frac{\mathbb{E}_2}{\Delta t^{\beta_2} \Gamma(2-\beta_2)},\quad C^M_2 = \frac{\mathbb{E}_2/\mathbb{E}_1}{\Delta t^{\beta_2 - \beta_1} \Gamma(2-\beta_2 + \beta_1)}
\end{equation}

\noindent\emph{Fractional Kelvin-Zener:}
\begin{align}
    C^{KZ}_1 = \frac{\mathbb{E}_2}{\Delta t^{\beta_2} \Gamma(2-\beta_2)},\quad & C^{KZ}_2 = \frac{\mathbb{E}_3}{\Delta t^{\beta_3} \Gamma(2-\beta_3)} \\
    C^{KZ}_3 = \frac{\mathbb{E}_2 \mathbb{E}_3/\mathbb{E}_1}{\Delta t^{\beta_2+\beta_3-\beta_1} \Gamma(2-\beta_1-\beta_3+\beta_2)},\quad & C^{KZ}_4 = \frac{\mathbb{E}_2/\mathbb{E}_1}{\Delta t^{\beta_2-\beta_1} \Gamma(2-\beta_2+\beta_1)}
\end{align}

\noindent\emph{Fractional Poynting-Thomson:}
\begin{align}
    C^{PT}_1 = \frac{\mathbb{E}_1}{\Delta t^{\beta_1} \Gamma(2-\beta_1)},\quad & C^{PT}_2 = \frac{\mathbb{E}_2}{\Delta t^{\beta_2} \Gamma(2-\beta_2)} \\
    C^{PT}_3 = \frac{\mathbb{E}_1 / \mathbb{E}_3}{\Delta t^{\beta_1-\beta_3} \Gamma(2-\beta_1+\beta_3)},\quad & C^{PT}_4 = \frac{\mathbb{E}_2/\mathbb{E}_3}{\Delta t^{\beta_2-\beta_3} \Gamma(2-\beta_2+\beta_3)}
\end{align}

\noindent\emph{Fractional Quasi-Linear viscoelastic:}
\begin{equation}
    C^{QLV}_1 = \frac{E A B}{\Delta t^{\beta} \Gamma(2-\beta)}
\end{equation}

\section{Return-Mapping Derivation for the Fractional Kelvin-Zener Model} \label{Ap:RMFKZ}

Recalling the discretized FKV model \eqref{eq:FKZdiscrete} employed as the viscoelastic part of the visco-elasto-plastic model:
\begin{dmath}\label{eq:FKVderiv_1}
    \sigma_{n+1} = (1+C^{KZ}_4)^{-1} \left[ C^{KZ}_1 \left(  \Delta\varepsilon^{ve}_{n+1}+\mathcal{H}^{\beta_2}\varepsilon\right)+C^{KZ}_2\left(\Delta\varepsilon^{ve}_{n+1}+\mathcal{H}^{\beta_3}\varepsilon^{ve} \right) + C^{KZ}_3\left(\Delta\varepsilon^{ve}_{n+1}+\mathcal{H}^{\beta_2+\beta_3-\beta_1}\varepsilon^{ve}\right)+C^{KZ}_4\left(\sigma_n-\mathcal{H}^{\beta_2-\beta_1}\sigma\right) \right],
\end{dmath}
where with the kinematic relationship \eqref{eq:discretizedkinematics} and the viscoplastic strain evolution \eqref{eq:discretized_vp}, we note that $\Delta\varepsilon^{ve}_{n+1} = \Delta\varepsilon_{n+1} - \Delta\gamma_{n+1}\sign(\sigma_{n+1})$. Therefore, \eqref{eq:FKVderiv_1} becomes:
\begin{dmath}\label{eq:FKVderiv_1_2}
    \sigma_{n+1} = (1+C^{KZ}_4)^{-1} \left[ C^{KZ}_1 \left(  \Delta\varepsilon_{n+1} +\mathcal{H}^{\beta_2}\varepsilon\right)+C^{KZ}_2\left(\Delta\varepsilon_{n+1} +\mathcal{H}^{\beta_3}\varepsilon^{ve} \right) + C^{KZ}_3\left(\Delta\varepsilon_{n+1}+\mathcal{H}^{\beta_2+\beta_3-\beta_1}\varepsilon^{ve}\right)+C^{KZ}_4\left(\sigma_n-\mathcal{H}^{\beta_2-\beta_1}\sigma\right) - \left(C^{KZ}_1+C^{KZ}_2+C^{KZ}_3\right)\Delta\gamma_{n+1}\sign(\sigma_{n+1}) \right].
\end{dmath}

Recalling the trial state for the FKZ model:
\begin{dmath}
    \sigma^{trial}_{n+1} = (1+C^{KZ}_4)^{-1} \left[ C^{KZ}_1 \left( \Delta\varepsilon_{n+1}+\mathcal{H}^{\beta_2}(\varepsilon-\varepsilon^{vp})\right)+C^{KZ}_2\left(\Delta\varepsilon_{n+1}+\mathcal{H}^{\beta_3}(\varepsilon-\varepsilon^{vp}) \right) + C^{KZ}_3\left(\Delta\varepsilon_{n+1}+\mathcal{H}^{\beta_2+\beta_3-\beta_1}(\varepsilon-\varepsilon^{vp})\right)+C^{KZ}_4\left(\sigma_n-\mathcal{H}^{\beta_2-\beta_1}\sigma\right) \right]
\end{dmath}
which, by combining the above two equations, we find:
\begin{dmath}\label{eq:FKVderiv_2}
\sigma_{n+1} = \sigma^{trial}_{n+1} - \sign(\sigma_{n+1}) \left(\frac{C^{KZ}_1+C^{KZ}_2+C^{KZ}_3}{1+C^{KZ}_4}\right) \Delta\gamma_{n+1}.
\end{dmath}

Finally, we obtain the loading/unloading sign consistency by following standard plasticity procedures:
\begin{equation}
    \sign(\sigma_{n+1})|\sigma_{n+1}| = \sign(\sigma^{trial}_{n+1})|\sigma^{trial}_{n+1}| - \sign(\sigma_{n+1}) \left(\frac{C^{KZ}_1+C^{KZ}_2+C^{KZ}_3}{1+C^{KZ}_4}\right) \Delta\gamma_{n+1},
\end{equation}
threfore,
\begin{equation}
    \sign(\sigma_{n+1})\left[|\sigma_{n+1}| + \left(\frac{C^{KZ}_1+C^{KZ}_2+C^{KZ}_3}{1+C^{KZ}_4}\right) \Delta\gamma_{n+1}\right] = \sign(\sigma^{trial}_{n+1})|\sigma^{trial}_{n+1}|,
\end{equation}
since both terms multiplying the sign functions on the left and right sides of the above equation are positive, we therefore conclude that $\sign(\sigma_{n+1}) = \sign(\sigma^{trial}_{n+1})$, and hence \eqref{eq:FKVderiv_2} becomes:
\begin{dmath}
\sigma_{n+1} = \sigma^{trial}_{n+1} - \sign(\sigma^{trial}_{n+1}) \left(\frac{C^{KZ}_1+C^{KZ}_2+C^{KZ}_3}{1+C^{KZ}_4}\right) \Delta\gamma_{n+1},
\end{dmath}
which completes the derivation.

\section*{Acknowledgments}
This work was supported by the ARO YIP Award (W911NF-19-1-0444), the NSF Award (DMS-1923201), the MURI/ARO Award (W911NF-15-1-0562), the AFOSR YIP Award (FA9550-17-1-0150), and NIH NIDCD K01DC017751 and R21DC020003.

\bibliographystyle{siamplain}
\bibliography{references}

\end{document}